\documentclass[11pt,3p,times,reqno]{amsart}
\usepackage[right=2.3cm,left=2.3cm,top=2.5cm,bottom=2.5cm]{geometry}
\usepackage[utf8]{inputenc}
\usepackage{amsmath,amsxtra,amssymb,latexsym,amscd,amsthm,mathrsfs,amsfonts,amstext,amsbsy,upgreek,dsfont,eufrak,float,cite}
\usepackage{fancyhdr,booktabs,verbatim}
\usepackage{multicol,graphicx,array,multirow,fancyhdr,booktabs,xcolor}
\usepackage{hyperref}
\hypersetup{colorlinks,linkcolor={red},citecolor={blue},urlcolor={cyan}}
\usepackage{bbm,csquotes}
\usepackage{extarrows}
\usepackage{esint,enumerate}
\usepackage{slashed}

\usepackage{mathtools,thmtools,pgf,tikz,graphicx}
\usepackage{cleveref}
\usepackage[shortlabels]{enumitem}
\declaretheorem[numberwithin=section]{theorem}
\usetikzlibrary{arrows}

\DeclareMathOperator*{\esssup}{{ess\,sup}} 
\def\nn{\nonumber}

\theoremstyle{definition} \newtheorem{definition}[theorem]{Definition}

\theoremstyle{plain} \newtheorem{lemma}[theorem]{Lemma}
\newtheorem{proposition}{Proposition}[section]

\newtheorem{remark}{Remark}[section]

\newcommand{\lf}{\left}
\newcommand{\rt}{\right}
\newcommand{\al}{\alpha}

\newcommand{\la}{\lambda}

\newcommand{\ep}{\epsilon}

\title[Long-time asymptotics and the modified conservation law for super-critical NLS]{\textbf{Long-time asymptotics of solutions and the modified  pseudo-conformal conservation law for super-critical nonlinear Schr\"odinger equation}}

\author[V.V. Au]{Vo Van Au}
\address[V.V. Au]{Division of Applied Mathematics, Science and Technology Advanced Institute, Van Lang University, Ho Chi Minh City, Vietnam; and Faculty of Applied Technology, School of Engineering and Technology, Van Lang University, Ho Chi Minh City, Vietnam}
\email{\url{vovanau@vlu.edu.vn}}

\author[T. Caraballo]{Tom\'{a}s Caraballo} 
\address[T. Caraballo]{Departamento de Ecuaciones Diferenciales y An\'{a}lisis Num\'{e}rico, Facultad de Matem\'{a}ticas, Universidad de Sevilla, C/Tarfia s/n, 41012-Sevilla, Spain}
\email{\url{caraball@us.es}}

\author[N.H. Tuan]{Nguyen Huy Tuan$^{*}$}
\address[N.H. Tuan]{Division of Applied Mathematics, Science and Technology Advanced Institute, Van Lang University, Ho Chi Minh City, Vietnam; and Faculty of Applied Technology, School of Engineering and Technology, Van Lang University, Ho Chi Minh City, Vietnam}
\email{\url{nguyenhuytuan@vlu.edu.vn}}

\thanks{$^*$Corresponding author: \url{nguyenhuytuan@vlu.edu.vn} (N.H. Tuan)}

\allowdisplaybreaks
\begin{document}
\date{\today}
\mdseries 
\maketitle
\begin{abstract}
In this paper, we discuss a class of nonlinear Schr\"odinger equations with the power-type nonlinearity: $(\mathrm i \frac{\partial}{\partial t}  + 
\Delta ) \psi =  \la |\psi|^{2\eta}\psi$ in $\mathbf R^N \times \mathbf R^+$. Based on the Gagliardo-Nirenberg interpolation inequality, we prove the local existence and long-time behavior (continuation, finite-time blow-up or global existence, continuous dependence) of the solutions to the $(H^q)$ super-critical Schr\"odinger equation. The corresponding scaling invariant space is homogeneous Sobolev $\dot H^{q_{crit}}$
with $q_{crit} > q$. Based on the estimates of the quadratic terms containing the phase derivatives used in the paper by Killip, Murphy and Visan \cite[SIAM J. Math. Anal. 50(3) (2018), 2681--2739]{KMV018} we shall study the stability with a stronger bound on the solutions to our problem. Moreover, from the arguments on virial-types presented in the paper by Killip and Visan \cite[Amer. J. Math. 132(2) (2010), 361--424]{KV010}, a modified pseudo-conformal conservation law is proposed. The Morawetz estimate for the solutions to the problem are also presented. 
		
		
		
	\end{abstract}
	\tableofcontents
	
	\section {Introduction and statement of the main results}
	The nonlinear Schr\"odinger (NLS for short) equation is used to describe a wavefield that is more complex than the normal form because its oscillation has a frequency proportional to the difference of the function value and its mean, see \cite{S26} and references therein. A solution of the Schr\"odinger equation that can be physically visualized is an amplitude function of a quantum mechanical particle that moves freely in $\mathbf R^N $ (see \cite{AM013,S26,KOPV017}). 
	In this study, we consider the following NLS equation with smooth power-type nonlinearity: 
	\begin{align}\label{Pro}
		\begin{cases}
			\lf(\mathrm i \dfrac{\partial}{\partial t}  + 
			\Delta \rt) \psi =  \la |\psi|^{2\eta}\psi, \quad \mbox{for}~~  x \in \mathbf R^N, ~ t \in \mathbf R^+,  \\
			\psi = \psi_0( x), \quad \mbox{for}~~ t=0,~  x \in \mathbf R^N, 
		\end{cases} \tag{$NLS$}
	\end{align}  
	where $N\geq 1$ and $\psi(x,t),~(x,t)\in \mathbf R^N \times \mathbf R^+$ is a complex-valued function; $\mathrm i = \sqrt{-1}$; 
	$|\psi(\cdot, t)|^2$ is the spatial probability density function; $\Delta = \partial^2/\partial x_1^2 + \cdots +\partial^2/\partial x_N^2$. The power-type nonlinearity $\la |\psi|^{2\eta} \psi,$ for $\la \in \mathbf R$ (or, possibly, $\la \in \mathbf C$) serves as a scaling factor, when $\la<0$ (or $\la>0$) the equation be called focusing (or defocusing); $1 \leq \eta\in \mathbf N$; the initial data $\psi_0$ belongs to $H^q(\mathbf R^N)$ for $\mathbf N \ni q > \frac{N}{2}$ and $q \geq \eta$.
	
With the power-type nonlinearity case $|\psi|^{2\eta}\psi ~ (\eta \geq 1)$, equation \eqref{Pro} becomes
scale invariant, that is, if $\psi(x, t)$ satisfies \eqref{Pro}, then so we have 
$
	\psi_\theta(x, t) = \theta^{\frac{1}{\eta}} \psi(\theta x, \theta^2t)$, with $\theta > 0.
$
For the critical exponent $q_{crit} = \frac{N\eta-2}{2\eta},$ and $\lf\|\psi_\theta(\cdot,t)\rt\|_{H^{q_{crit}}(\mathbf R^N)} = \lf\|\psi(\cdot,\theta^2 t)\rt\|_{H^{q_{crit}}(\mathbf R^N )}$,
if we consider the case of $q_{crit} > q \Longrightarrow \eta > \frac{2}{N-2q}$ and $\psi_0 \in H^q(\mathbf R^N)$ then the Problem \eqref{Pro} is called the $(H^q)$ super-critical problem. Moreover, if $q_{crit} >0$ then the problem is called $(L^2)$/mass super-critical case, if $q_{crit}<1$ then the problem is called $(H^1)$/energy sub-critical case (see \cite{MR008,CTG009}). Problem \eqref{Pro} with condition $q>\frac{N}{2}$ called the $(H^q)$ super-critical NLS equation.	
	
	It is well known that  smooth solutions to the NLS equation \eqref{Pro} satisfy the following conservation laws for all $t > 0$ (see \cite{MR008,CGT010})
\begin{align} \label{ma}
	L^2-\mbox{norm}: \quad \mathbb T[\psi](t) = \int_{\mathbf R^N} |\psi (x,t)|^2 dx = \int_{\mathbf R^N} |\psi_0(x)|^2 dx = \mathbb T[\psi_0]; \tag{Mass} 
\end{align}
\begin{align} \label{ener}
	\mathbb E[\psi](t) &= \frac{1}{2} \int_{\mathbf R^N}  |\vec \nabla \psi (x,t) |^2 dx + \frac{\la}{2\eta+2} \int_{\mathbf R^N} |\psi (x,t)|^{2\eta+2}  dx \nn\\
	&= \frac{1}{2} \int_{\mathbf R^N}  |\vec \nabla \psi_0 (x) |^2 dx + \frac{\la}{2\eta+2} \int_{\mathbf R^N} |\psi_0 (x)|^{2\eta+2}  dx = \mathbb E[\psi_0]; \tag{Energy} 
\end{align}
\begin{align} \label{mom}
	\mathbb M[\psi](t) = \Im \int_{\mathbf R^N} \bar \psi (x,t) \vec \nabla \psi (x,t) dx = \Im \int_{\mathbf R^N}  \bar \psi_0(x) \vec \nabla \psi_0(x) dx = \mathbb M[\psi_0]. \tag{Momentum} 
\end{align}
Here, $\Im(z)$ is the imaginary part of the complex number $z$ and $\bar{\psi}$ is the complex conjugate term of $\psi$. More precisely, $\varPsi= |\psi|^2$ (mass density) and the \textit{local} conservation of mass reads (see \cite[page 12]{CTG009})
\begin{align} \label{mass}
	\mathbb T_{\mathrm{loc}}[\psi](x,t) &= \frac{\partial \varPsi }{\partial t}(x,t)  + 2 \mbox{div} \lf(\overrightarrow{Q^\nabla}(\psi)(x,t) \rt) \nn\\
	&= \frac{\partial \varPsi }{\partial t}(x,t)  + 2 \nabla_1 Q^\nabla_1(\psi)(x,t) + 2 \nabla_2 Q^\nabla_2(\psi)(x,t) \nn\\
	&\quad \stackrel{\cdots}{+} 2 \nabla_{N-1} Q^\nabla_{N-1}(\psi)(x,t) + 2\nabla_N Q^\nabla_N (\psi)(x,t) = 0.   \tag{$\mbox{loc-Mass}$}
\end{align}
Here, $\varPsi  = |\psi  |^2 = \psi  \bar{\psi}$ and we note that 
	\begin{align} \label{def:Q}
	\overrightarrow{Q^\nabla} (\psi) = \lf(Q^\nabla_1 (\psi), ...,Q^\nabla_N(\psi)\rt) \quad \mbox{and}\quad Q^\nabla_m(\psi)  = \Im \lf(\bar{\psi} \nabla_m \psi \rt),~ m = 1,...,N.	
	\end{align}
The \textit{local} momentum conservation is (see \cite[page 12]{CTG009}) 
\begin{align} \label{Mom}
	\mathbb M_{\mathrm{loc}}[\psi](x,t)
	= \frac{\partial Q^\nabla_m (\psi)}{\partial t} (x,t) &+  \nabla_1 \lf(\mathbf 1_{m, 1} L^\nabla(\psi)  + S^\nabla_{m, 1} (\psi) \rt)(x,t) \nn\\
	& + \nabla_2 \lf(\mathbf 1_{m, 2} L^\nabla(\psi) + S^\nabla_{m, 2} (\psi) \rt)(x,t)   \nn\\
	& \stackrel{\cdots}{+} \nabla_{N-1} \lf(\mathbf 1_{m, N-1} L^\nabla (\psi)  + S^\nabla_{m, N-1} (\psi)\rt)(x,t) \nn\\
	& + \nabla_N \lf(\mathbf 1_{m, N} L^\nabla(\psi) + S^\nabla_{m, N} (\psi) \rt)(x,t)  =0,  \tag{$\mbox{loc-Momentum}$}
\end{align} 
where $m = 1,...N$ and $\mathbf 1_{m,n}$ is the Kronecker delta function 
\begin{align} \label{Kronecker}
	\begin{cases}
		\mathbf 1_{m,n}  = 1, \quad \mbox{if}\quad m=n; \\
		\mathbf 1_{m, n} = 0, \quad \mbox{if} \quad m \neq n,
	\end{cases} \quad 1 \leq m,n \leq N,
\end{align}
and the Lagrangian density is
\begin{align} \label{lar}
	L^\nabla(\psi)  = -\frac{1}{2} \mbox{div}\lf(\vec \nabla |\psi|^2 \rt) + \frac{\la \eta}{\eta+1} \psi^{2\eta+2} ,
\end{align} 
where $\vec \nabla = (\nabla_1,...,\nabla_N)$. The symmetric tensor is
\begin{align} \label{ten}
	S^\nabla_{m,n}(\psi)  = 2  \Re \lf(\nabla_m \psi   \nabla_n \bar{\psi} \rt), \quad \mbox{for}~ m = \overline{1,N}, n = \overline{1,N},
\end{align}
with $\Re(z)$ be the real part of the complex number $z$. 	

The pseudo-conformal conservation (PCC) law was first discovered by Ginibre and Velo in the series of papers \cite{GV79-I,GV79-II,GV80}, considered as an essential tool to study the asymptotic behavior of solutions for NLS equations.
J. Bourgain \cite{Bo99} considered the $3D$ defocusing NLS with $\eta =2$
	\begin{align} \label{Bo99}
	\begin{cases}
	\lf(\mathrm i \dfrac{\partial}{\partial t}  + 
	\Delta \rt) \psi - |\psi|^{4}\psi = 0, \quad  x \in \mathbf R^3, ~ t \geq t_0,  \\
	\psi|_{t=t_0} = \zeta \psi(t_0), \quad x \in \mathbf R^3,
	\end{cases}
	\end{align}
where $0 < \zeta < 1$ is a radial bump function chosen and the author proved that the solutions of problem satisfy the following PCC law   \begin{align} \label{pse-Bo}
\lf\|[(\cdot) +2\mathrm i(t-t_0)\vec \nabla] \psi (\cdot, t) \rt\|_{L^2(\mathbf R^3)}^2& + \frac{4 }{3} (t-t_0)^2 \lf\|\psi(\cdot,t)\rt\|_{L^6(\mathbf R^3)}^6\nn\\
 &= \lf\|(\cdot)\psi(\cdot,t_0)\rt\|_{L^2(\mathbf R^N)}^2 - \frac{16}{3}
\int_{t_0}^t s \int_{\mathbf R^3} |\psi(x,s)|^6 dxds \nn\\ 
&\leq \lf\|(\cdot) \psi (\cdot,t_0)\rt\|_{L^2(\mathbf R^3)}^2.
\end{align} 
By considering Problem \eqref{Bo99}, J. Colliander et al. also pointed out PCC \eqref{pse-Bo} in paper \cite[Eq. (8.15)]{CKSTT008}.
Obviously, Problem \eqref{Pro} is more general than Problem \eqref{Bo99}, then PCC law for this problem is much more difficult to establish. As far as we know, in general, the PCC law for Problem \eqref{Pro} has not been studied extensively. There are two formal ways to set up the PCC for Problem \eqref{Pro}. One is using the pseudo-conformal symmetry, indeed, for $\eta = \frac{2}{N}$ and if $\psi$ is a solution of \eqref{Pro} then
	$\varTheta(x,t)=\frac{1}{|t|^{\frac{N}{2}}} \overline{\psi\lf(\frac{x}{t},\frac{1}{t}\rt)} e^{\frac{i|x|^2}{4t}}$
is also a solution to \eqref{Pro} for $t \neq 0$. The function $\varTheta$ satisfies the equation
	$$\mathrm i\frac{\partial \varTheta}{\partial t}  + \Delta \varTheta = \la t^{\eta - 2} |\varTheta|^{2\eta} \varTheta,$$
and by direct arguments on this equation, we obtain the PCC for Problem \eqref{Pro}, see \cite[Subsection 7.5, page 224]{C003}. Another way is to use the virial identities and this is also our aim to prove the following quantity 
\begin{align} \label{def-Pseudo}
	\forall t \geq 0:	~~ \mathbb P[\psi](t)  &= \|[(\cdot) + 2\mathrm it \vec\nabla] \psi (\cdot,t)\|_{L^2 (\mathbf R^N)}^2  \nn\\
	& - \frac{4t^2\la N}{N+2}\int_{\mathbf R^N}  |\psi(x,t)|^{\frac{4}{N}+2}  dx  - 2N \int_{\mathbf R^N} \mathrm{div} \lf(\vec \nabla \mathcal H(\psi)(x,t)  \rt) dx ,
\end{align}
is conserved. Here, the function $\mathcal H$ satisfies some necessary conditions:
\begin{align} \label{cond:H}
	\begin{cases}
		\mathcal H(\psi)(x,0) = 0; \\  
		\dfrac{\partial \mathcal H(\psi)}{\partial t} (x, t) = t |\psi|^2, 
	\end{cases} \mbox{for a.e. on}~ \mathbf R^N.
\end{align}	
Based on the properties of the virial identity which is called $(1-1)$ Morawezt action (see \autoref{virial} for its definition), we prove a the modified by considering the new version of the PCC law to Problem \eqref{Pro} as shown in \eqref{def-Pseudo} in $\mathbf R^N$, i.e.,
\begin{align} \label{PCC}
\mathbb P[\psi](t) ~~ \mbox{is conserved} \quad \Longleftrightarrow \quad \mathbb P[\psi](t) = \mathbb P[\psi_0], \quad \forall t>0. \tag{\mbox{PC}}
\end{align}
\begin{remark}
Since Problem \eqref{Pro} is more general than \eqref{Bo99}, it can be compared directly to see that \eqref{PCC} includes \eqref{pse-Bo}. The conditions in \eqref{cond:H} of $\mathcal H(\psi)$ are satisfied for a class of functions containing the solution of the problem under some time-smooth conditions. This is an extension of the previous result to show that there are more nonlinear functions containing particles that still guarantee conservation. However, these functions must satisfy the constraint conditions as \eqref{cond:H}. 
\end{remark} 

The main  equation in \eqref{Pro} is a particular case of the general NLS equation
	\begin{align} \label{NLS1}
		\mathrm i \dfrac{\partial \psi}{\partial t} = - \Delta \psi +  F (\psi).
	\end{align}
	Equation \eqref{NLS1} has been studied by V. Banica \cite{B004} for $F (\psi) = - |\psi|^{\eta-1}\psi$, the author was concerned with the focusing NLS equation posed on a regular domain $\Omega$ of $\mathbf R^N$, with Dirichlet boundary conditions. The blow-up rates are bounded. Moreover, they proved that if the blow-up occurs on the boundary, the blow-up rate grows faster than $(T - t)^{-1}$. For dimension $N \geq 3$, $F (\psi) = - |\psi|^{\eta-1}\psi$ for $1 + \frac{N}{4}< \eta < \frac{N+2}{N-2}$,  F. Merle and P. Rapha\"{e}l \cite{MR008} proved the blow up of the critical norm for some radial $L^2$ super-critical NLS equations.
	
	In dimension $2$, Y. Wang \cite{W008} discussed the NLS equation containing the general potential
	$F (\psi) =  V(x) |\psi|^{\eta-1}\psi$
	with some assumptions of potential $V(x)$. The paper obtained some sharp conditions for global existence and blow-up of solutions in $H^1(\mathbf R^2)$. For the focusing cubic NLS, C. Sulem et al. \cite{CRSW004} proved the intermediate case of blow-up solutions from initial data $\psi_0 \in H^q(\mathbf R^N)$, with a very strict index satisfying $1 > q > q_s$, for $q_s \leq \frac{1}{5} (1 + \sqrt{11})$, it called below $H^{1}$ cases. We refer the readers to the excellent works of C. Sulem et al. in \cite{SS97,BS003,OS012}.
	
L. Vega et al. \cite{KPPV020,EKPV010,VB011} studied the unique continuation of solutions to non-local non-linear dispersive equations and related to NLS. In particular, Banica and Vega \cite{BV020-1,BV020-2} demonstrated the existence of a unique solution of the binormal flow with datum a polygonal line. Some other outstanding works of L. Vega and his group on dispersive equations can be suggested to the reader as \cite{AV019,GV013,AMV013} and some related references therein. Many researchers have studied  the asymptotic behavior of  solutions to NLS equations: G.M. Bisci and V.D. R\u{a}dulescu \cite{BR015}, N. Burq et al. \cite{BGT003}, B.J. Campos and P.I. Naumkin \cite{CN021}, T. Oh et al. \cite{OW020}, H. Xu \cite{Xu017}, Q. Ding et al. \cite{DW018}, A.J. Fern\'{a}ndez et al. \cite{FW021}, L. Jeanjean et al. \cite{JL021}, R. Killip et al. \cite{KMVZZ018}, T. Caraballo et al. \cite{ZCL021}, K. Pravda-Starov \cite{MS021}, D.V. Duong \cite{D018,D019} and the references therein.  To our knowledge, there are many works on the $(L^2)$/mass super-critical \cite{MR008,B004,CRSW004} and $(H^1)$/energy sub-critical \cite{CRSW004,KV010,OT91} problems, but the results on long-time asymptotics for super-critical NLS problems are still limited. Our results on $(H^q)$ super-critical include the results on energy super-critical problems.

In this paper, to solve Problem \eqref{Pro}, we have tried to overcome the following challenges: \\
$\bullet$ First, without the help of Strichartz estimates, but using the Gagliardo-Nirenberg (G-N) interpolation  inequality (see \autoref{GN} below) we obtain the local well-posedness results of the solutions: local existence, continuation, global existence (or finite time blow-up) and continuous dependence on initial data. In \cite{MR008}, the author used Strichartz estimates to obtain the finite-time blow-up and lower bound for the mass super-critical case in  Sobolev norm. This is the reason that motivates us to study the long-term behavior for the $(H^q)$ super-critical of \eqref{Pro}. However, we do not rely on the usefulness estimates of Strichartz used in \cite{MR008} and many other works. Thus, we face the difficulty of establishing the Lipschitz property of the nonlinear source function $\la |u|^{2\eta}u$ to be able to use a complex tool as the G-N interpolation inequality. It is worth mentioning here that we need the conditions of the power $\eta$ satisfy $\mathbf N \ni \eta \in [1,q]$.  It is a necessary condition to establish techniques for the application of G-N interpolation inequality. Let us cover a little bit of the technique that is closely related to use source functions of the form $\mathcal F_{\eta}(\psi) = \la |\psi|^{2\eta}\psi$. We need to estimate the following term $\|\mathcal F_{\eta}(\psi)\|_{H^q(\mathbf R^N)}$ and   
	\begin{align*}
	\|\mathcal F_{\eta}(\psi)\|_{H^q(\mathbf R^N)} = \lf(\sum_{|\beta| \leq q} \sum_{\eta \leq q} \Big\|\mathcal F_{\eta}^{(\eta)}(\psi) \prod_{j=1}^{\eta} \mathscr D^{\al_j} (\psi) \Big\|_{L^2(\mathbf R^N)}^2 \rt)^{\frac{1}{2}}.
	\end{align*}
From the above formula, we have the following two important problems. First, the estimation of the norms of the term $\mathcal F_{\eta}^{(\eta)}(\psi)$ explains why we need to consider the term $\mathcal F_{\eta}(\psi) = \la |\psi|^{2\eta}\psi$ (i.e., related to the exponent $2\eta$ of $|\psi|$) instead of other nonlinear functions. It is a clear advantage that we can use the hypothesis \eqref{H2}. Second, estimating the norms of the term $\prod_{j=1}^{\eta} \mathscr D^{\al_j} (\psi)$ requires the help of G-N interpolation inequality, and it is interesting to choose the complex coefficients of this inequality related to $\eta$ and $q$. It also has to satisfy the conditions of H\"older inequality and generalized H\"older inequality. Through Sobolev embeddings, we establish the properties of the solution on a Sobolev space $H^q(\mathbf R^N)$ with the conditions $\frac{N}{2} < q \in \mathbf N$ and $q\geq \eta$.  These results are presented in \autoref{existence} and \autoref{blow-up}, below. \\ 
$\bullet$ The second challenge is to estimate a stronger bound space-time for solutions. In \cite{KMV018}, the authors estimate the ``phase derivative" terms (see \cite[in Strategy of the proof]{KMV018}), such terms are significantly more difficult to estimate. These estimates contain higher order gradients with $\widehat{|\nabla|^sf}(\xi) = |\xi|^s \widehat f(\xi)$. Also, based on ideas in the paper \cite{KMV018}, however, we do not use the estimations for phase derivative terms, but instead, we consider the expansions of generalized virial quantity with two particles, i.e. contains both $\psi(x,t)$ and $\psi(y,t)$ (see \cite[Eq. (1.8)]{CTG009})
	\begin{align} \label{Mora2}
	\mathcal V_\delta(t) = \iint_{\mathbf R^N \times \mathbf R^N}  |\psi(y,t) |^2  \vec{\slashed{\nabla}}_x \delta(x-y) \cdot \overrightarrow{Q^\nabla} (\psi)(x,t) dx dy,
	\end{align}
where $\delta$ is a weighted function with two variables $(x,y) \in \mathbf R^N \times \mathbf R^N$. By choosing $\delta(x-y) = |x-y|$, $x,y \in \mathbf R^N$, the estimations of this virial \eqref{Mora2} are closely related to the conservation of mass \eqref{ma} and the local conservation  laws \eqref{mass}, \eqref{Mom}. This result is presented in \autoref{stability} and \autoref{prop:space-time-est}.
\begin{remark} \label{virial} Some notes on the virial identity are as follows.
	\begin{itemize}
\item[(i)] The generalized virial quantity \eqref{Mora2} can be called the two-particles Morawetz action with two-variables weighted function.
	We rewrite \eqref{Mora2} as follows
	\begin{align*}
		\mathcal V_\delta(t) = \int_{\mathbf R^N}  |\psi(y,t) |^2  \mathcal M_\delta(y,t) dy, \quad \mbox{with}~ \mathcal M_\delta(y,t) = \int_{\mathbf R^N}  \vec{\slashed{\nabla}}_x \delta(x-y) \cdot \overrightarrow{Q^\nabla} (\psi)(x,t) dx .
	\end{align*}
	By using a reduced form of \eqref{Mora2} for the one-variable weighted function $\delta :=\delta(x)$, $x\in \mathbf R^N$, the restrict virial 
	\begin{align} \label{Mora1}
		\mathcal M_\delta(t) \equiv \mathcal M_\delta(0,t) = \int_{\mathbf R^N}  \vec{\slashed{\nabla}}_x \delta(x) \cdot \overrightarrow{Q^\nabla} (\psi)(x,t) dx,
	\end{align}
	can be called the one-particle Morawetz action with one-variable weighted function (see \cite[Eq. (1.6)]{CGT010}). 
\item[(ii)] From now on, we shall say the generalized virial $\mathcal V_\delta(t)$ in \eqref{Mora2} is $(2-2)$ Morawetz action (i.e., two-particles and two-variables weighted function) and the restrict virial $\mathcal M_\delta(t)$ in \eqref{Mora1} is $(1-1)$ Morawetz action (i.e., one-particle and one-variable weighted function).
\item[(iii)] In some situations, another form of virial quantity was also considered in \cite[page 29]{KV010} with 
	$$\int_{\mathbf R^N} \delta(x) |\psi(x,t)|^2dx, \quad\mbox{with the weighted function}~~ \delta(x) = \phi (\frac{|x|}{R}), \quad \phi(r) = \begin{cases}
		1, \quad \mbox{for}~ r\leq 1,\\
		0, \quad \mbox{for}~ r \geq 2.
	\end{cases}	$$
\item[(iv)] The minimal condition to impose on the weighted function $\delta :=\delta(x-y)$ in \eqref{Mora2} is that the matrices of second partial derivatives 
$$\lf(\begin{array}{cccc}
\slashed\nabla_{x,1}^2 \delta& \slashed\nabla_{x,1}\slashed\nabla_{x,2} \delta &\dots&\slashed\nabla_{x,1}\slashed\nabla_{x,N} \delta  \\
\slashed\nabla_{x,2}\slashed\nabla_{x,1} \delta& \slashed\nabla_{x,2}^2 \delta &\dots&\slashed\nabla_{x,2}\slashed\nabla_{x,N} \delta\\
\vdots & \vdots&\ddots&\vdots\\
\slashed\nabla_{x,N}\slashed\nabla_{x,1} \delta&\slashed\nabla_{x,N}\slashed\nabla_{x,2} \delta&\dots&\slashed\nabla_{x,N}^2 \delta
\end{array}	\rt)_{N\times N}$$
and 
$$\lf(\begin{array}{cccc}
\slashed\nabla_{y,1}^2 \delta& \slashed\nabla_{y,1}\slashed\nabla_{y,2} \delta&\dots&\slashed\nabla_{y,1}\slashed\nabla_{y,N} \delta \\
\slashed\nabla_{y,2}\slashed\nabla_{y,1} \delta& \slashed\nabla_{y,2}^2 \delta&\dots&\slashed\nabla_{y,2}\slashed\nabla_{y,N} \delta\\
\vdots & \vdots&\ddots&\vdots\\
\slashed\nabla_{y,N}\slashed\nabla_{y,1} \delta&\slashed\nabla_{y,N}\slashed\nabla_{y,2} \delta&\dots&\slashed\nabla_{y,N}^2 \delta
\end{array}	\rt)_{N\times N}$$
are positive definite. With $\delta = |x-y|$, see Appendix \ref{A3} for the proof that the matrices $$\lf\{\slashed\nabla_{x,m}\slashed\nabla_{x,n} |x-y|\rt\}_{m,n=1}^N ~\mbox{and}~ \lf\{\slashed\nabla_{y,m}\slashed\nabla_{y,n} |x-y|\rt\}_{m,n=1}^N ~~\mbox{are positive definite}.$$ 
	\end{itemize}
\end{remark}
\noindent $\bullet$ The next challenges are to establish a new version (modified) of the PCC law related to Problem \eqref{Pro}. By choosing $\delta(x) = |x|^2$ in $(1-1)$ Morawetz action and combining with the conservation laws \eqref{ener}, \eqref{mass} and \eqref{Mom}, we prove the solutions of NLS \eqref{Pro} satisfy the modified PCC law \eqref{PCC}. To do this, we have combined the complex analytical properties of $\mathcal M_{\delta=|x|^2}(t)$ to obtain a consistent result in a conservative case that depends on the selection of power $\eta$ of the nonlinear source function related to $N$-dimensional real space (here, we pick $\eta N =2$ with $N=1$ or $N=2$). Moreover, in one dimension of \eqref{Pro}, as a consequence of the $(1-1)$ Morawetz action for $\eta \in [1,2)$, the Morawetz estimate and decay estimate on time are proposed. These results are provided in \autoref{pseudo-conformal}, \autoref{Morawetz estimate} and \autoref{decay}, below. 
	
	Now, we shall indicate the main results of this paper. The first theorem we want to mention is the result on the existence of a local solution, from which we prove that this problem can have solutions over longer time intervals (continuation of the solutions, see \autoref{def-cont} below). More specifically, we consider the following theorem:
	
	\begin{theorem}[$H^q$-local and large-time existence] \label{existence} For $N\geq 1$, $\mathbf N \ni \eta \geq 1$, there exists an integer number $q$ satisfying $q > \frac{N}{2}$ and $q \geq \eta$. Let $\psi_0 \in H^q(\mathbf R^N)$ and the nonlinearity $\la |\psi|^{2\eta}\psi$ satisfies hypotheses \eqref{H1}-\eqref{H2} below. Then:
		\begin{itemize}
			\item[(I)] There exists a positive constant $T$ (depending  only on $\psi_0$) such that  Problem \eqref{Pro} admits a (unique) mild solution in $L_t^\infty H^q_x (\mathbf B^N)$ for $\mathbf B^N := (0,T) \times \mathbf R^N$.
			\item[(II)] The unique solution of Problem \eqref{Pro} on $[0,T]$ can be stretched to the interval $[0,T_+]$ for some constants $T_+> T $, so the extension function is also the (unique) mild solution on $[0,T_+] $ of Problem \eqref{Pro}. 
		\end{itemize} 
	\end{theorem}
	
	\begin{remark} A more natural and interesting way to consider NLS is to avoid the dependence on Strichatrz estimates. For dimension $\geq 1$, we have tried to set up solution properties of Problem \eqref{Pro} on $H^q(\mathbf R^N)$ by using the Gagliardo-Nirenberg interpolation inequality. Although there is still a limit to the index of $q$ (namely, has to be an integer and $q>\frac{N}{2}$), we have achieved encouraging results.
	\end{remark}
	
	\begin{remark}
		For \autoref{existence} and \autoref{blow-up}, our results are limited to the cases $2q \leq N$ or $q$ is not an integer number, because of the following two difficulties:
		\begin{itemize}
			\item[i)] For $2q \leq N$, there does not exist the Sobolev embedding $H^q(\mathbf R^N) \hookrightarrow L^\infty(\mathbf R^N) $ (see \autoref{SE}).
			\item[ii)] For $q \notin \mathbf N$, we have more difficulties to estimate the Lipschitz property in $H^q(\mathbf R^N)$ of the reaction terms. 
		\end{itemize}
	\end{remark}

	As an immediate consequence of (II) in \autoref{existence}, we guarantee the existence of a
	maximal time (see \autoref{def-max}). Next, we prove the results on non-continuation due to a blow-up (see \autoref{def-blow}) or global existence. Furthermore, the continuous dependence on the initial data is also presented.

\begin{theorem}[Blow-up alternative, continuous dependence in $H^q$] \label{blow-up}
Assume that the conditions in \autoref{existence} hold. For every $\psi_0 \in H^q(\mathbf R^N)$ there exists a maximal time $T^{\star} > 0$ and a unique solution $\psi \in L_t^\infty H_x^q (\mathbf B^N_\star)$ with $\mathbf B^N_\star := (0,T_\star) \times \mathbf R^N$ which satisfies \eqref{exact-sol} and, in addition, satisfies the following conclusions:
	\begin{itemize}
	\item[(I)] If $T_\star < \infty$, then $\lim_{t \uparrow T_\star^- }
	\lf\|\psi(\cdot,t)\rt\|_{H^q (\mathbf R^N)} = \infty$, 
	or global existence if $T_\star = \infty$.
	\item[(II)] If $\psi_{0j} \to \psi_0$ in $H^q(\mathbf R^N)$ and if $\psi_{j}$ is the maximal solution corresponding to the initial data $\psi_{0j}$, then $\psi_j \to \psi$ in $L_t^\infty H_x^q(\mathbf B^N)$ for every interval $[0,T] \subset [0,T_\star)$.
	\end{itemize}
\end{theorem}
We state the following theorem to estimate the stability of solutions.	
	\begin{theorem}[Stability] \label{stability}
	For $N\geq1$, assume that Problem \eqref{Pro} has a solution $\psi \in C(\mathbf R^+,C_0^\infty (\mathbf R^N))$ with $\psi(x,0) = \psi_0 \in C_0^\infty (\mathbf R^N)$. Then the following estimate holds
	\begin{align} \label{est:space-time}
	\frac{C_{N,\pi}(N-1)}{2} \lf\|(- \Delta)^{\frac{3-N}{4}} |\psi|^2\rt\|_{L_t^2L_x^2(\mathbf R^+ \times \mathbf R^N)}^2 +  \frac{\la \eta(N-1)}{\eta+1} \iiint_{\mathbf R^+ \times \mathbf R^N \times \mathbf R^N} \frac{ |\psi(y,t)|^2 \lf| \psi (x,t)\rt|^{2\eta+2}  }{|x-y|}   dxdydt \nn\\
	\leq  \lf\|\psi_0 \rt\|_{L^2(\mathbf R^N)}^2 \sup_{t \in \mathbf R^+ \atop y \in \mathbf R^N} \lf|\int_{\mathbf R^N} \vec{\slashed{\nabla}}_x |x-y| \cdot \overrightarrow{Q^\nabla} (\psi)(x,t) dx \rt|.
	\end{align}
\end{theorem}	
The following proposition is a sharpening
of the results in \autoref{stability}.
\begin{proposition} \label{prop:space-time-est}
For dimension $\geq 1$, the solutions of the defocusing case ($\la>0$) of Problem \eqref{Pro} satisfy the following estimate
	$$\lf\|(- \Delta)^{\frac{3-N}{4}} |\psi|^2\rt\|_{L_t^2L_x^2(\mathbf R^+ \times \mathbf R^N)}  \lesssim \sqrt[4]{\mathbb T^{3}[\psi_0]} \sqrt{\mathbb E[\psi_0]}.$$  
\end{proposition}	
	A special case of the conservation law for Problem \eqref{Pro} involving the exponent of the reaction term with $\eta=2/N$ is considered. From the condition $\eta \geq 1$ in \autoref{existence} implies that we shall deal Problem \eqref{Pro} with $N=1$ or $N=2$. We present the following theorem.
	
	\begin{theorem}[Modified pseudo-conformal conservation law] \label{pseudo-conformal} For $N=1,2$, suppose that the conditions in \autoref{existence} hold.  Then  Problem \eqref{Pro} possesses a unique solution  $\psi \in L_t^\infty H_x^q(\mathbf B^N)$. For $\eta N= 2$, we have the pseudo-conformal conservation law
		\begin{align} \label{Pseudo}
			\mathbb P[\psi](t) = \mathbb P[\psi](0)=\|(\cdot)  \psi_0\|_{L^2 (\mathbf R^N)} , \quad t>0,  
		\end{align}
		where $\mathbb P[\psi](t)$ be defined as in \eqref{def-Pseudo}.
	\end{theorem}

	When $\eta N<2$ and combined with the condition $\eta \geq 1$ in \autoref{existence} then we infer that the Problem \eqref{Pro} is considered for $N=1$ and $\eta \in [1,2)$, we have space-time estimation of the solutions, which is shown in the following theorem.
	
	\begin{theorem}[Morawetz estimate] \label{Morawetz estimate} In one dimension of Problem \eqref{Pro}, assume that  conditions in \autoref{existence} are fulfilled. Then for $1 
		\leq \eta < 2$,  we have
		\begin{align} \label{Est}
			\iint_{[0,T]\times \mathbf R } t |\psi(x,t)|^{2\eta +2} dx dt \lesssim  \sup_{0 \leq t \leq T} \mathbb P[\psi](t).
		\end{align}
	\end{theorem}
	From \autoref{Morawetz estimate}, we investigate the following proposition.
	\begin{proposition}[Time decay estimate] \label{decay}
		For $N=1$, $\eta \in [1,2)$, suppose that $$\sup_{0< t \leq T} \lf( \frac{d}{dt} \mathbb P[\psi](t) \rt) \leq C,$$ for some $C>0$. Then, the decay estimate of $\psi$ in $L^{2\eta+2}(\mathbf R)$  holds for all $t \in (0,T]$, with
		$$\lf\|\psi(\cdot,t)\rt\|_{L^{2\eta+2} (\mathbf R)} \lesssim t^{-\frac{1}{2\eta+2}}, \quad t>0.$$
	\end{proposition}
	
	\section{Mathematical background and some related tools}
	First, we set up some basic notations and tools that are useful for our work.  For the quantities $f, g>0$, we shall write $f \lesssim g$ if there exists a constant $C > 0$
	such that $f \leq Cg$. The notation $\vec{\nabla}$ denotes the gradient of a scalar-valued differentiable function at the point $x \in \mathbf R^N$ and $\vec{\nabla}x = (\nabla_{1}x,...,\nabla_{N}x) = (\frac{\partial x}{\partial x_1},..., \frac{\partial x}{\partial x_N})$. In some cases, to discriminate the gradient of two spatial variables $x=(x_1,...,x_N)$ and $y$, we also write $\vec{\slashed \nabla}_x$ and  $\vec{\slashed \nabla}_y$.  In addition, $\stackrel{\cdots}{\pm}=\pm \cdots \pm$ and $\vdots \pm$ be the calculation $\pm$ vertical expansion. Let us set $[u_n]_{n=1}^N = (u_1,\cdots u_N)$ and $[\int \cdots \int][u_n]_{n=1}^N = \int u_1 \stackrel{\cdots}{+}
	\int u_N$.
	
	The notation $\|\cdot\|_{X}$ stands for the norm in the Banach space $X$. For $1  \leq m \leq \infty, ~ T >0,$ the Banach space of real-valued measurable functions $f: (0,T) \equiv I \to  X$ can be well defined with norms 
	\begin{align*}
		\|f\|_{L^m(I; X)} = \lf\{
		\begin{array}{ll}
			\lf( \displaystyle \int_{I} \lf\|f ( t) \rt\|^m_{X}   dt \rt)^{\frac{1}{m}} < \infty ,\quad  \mbox{for}~  1 \leq m < \infty,\\
			\displaystyle \esssup_{t\in I} \lf\|f ( t) \rt\|_{X} < \infty, \quad \mbox{for}~  m = \infty,
		\end{array} \rt. 
	\end{align*}
	From now on, assume that the spatial space $X$ has $N$-dimensions, we shall write $L^m_tX_x(\mathbf B^N)$ with $\mathbf B^N = I \times \mathbf R^N$ instead of the notation $L^m(I; X(\mathbf R^N))$. The norm of the function space $C^k(\bar I; X),$ for $0 \leq k \leq \infty$ is denoted by
	\begin{align}
		\|f\|_{C^k(\bar I; X)} =
		\sum_{j=0}^k
		\sup_{t \in \bar I}
		\lf\|f^{(j)} (t) \rt\|_{X} < \infty,
	\end{align}
	where $f^{(j)}$ is called the $j-$ th derivative of the function $f$. If $k = 0$ and $I$ is bounded interval then the space $C(I;X)$ is a Banach space with the norm of $L^\infty(I;X)$. For $k = \infty$, the space $C^\infty(\bar I;X)$ consists of functions with continuous strong derivatives of all orders.
	
	For a positive integer $q$ and $p \in [1,\infty)$, the Sobolev space $W^{q,p}(\mathbf R^N)$ consists of all locally integrable functions $f: \mathbf R^N \to \mathbf R$ such that the weak derivatives of the order less than or equal to $q$ belong to $L^p(\mathbf R^N)$, i.e.
	$\mathscr D^\beta (f) \in L^p(\mathbf R^N), ~ \mbox{for}~ 0 \leq |\beta| \leq q.$
	The Sobolev space $W^{q,p}(\mathbf R^N)$ is a Banach space when associated with the norm
	\begin{align*}
		\lf\|f \rt\|_{W^{q,p} } =\lf\{\begin{array}{ll}
			\displaystyle \lf(\sum_{|\beta| \leq q} \int_{\mathbf R^N}  |\mathscr D^\beta (f)  |^p dx\rt)^{\frac{1}{p}} < \infty, \quad \mbox{for}~ 1 \leq p < \infty,\\
			\displaystyle \max_{|\beta| \leq q} \sup_{x \in \mathbf R^N}   |\mathscr D^\beta (f)  | < \infty, \quad \mbox{for}~ p = \infty .
		\end{array} \rt.  
	\end{align*}
	Note that, for $p=2$, the space $H^q(\mathbf R^N)$ is a Hilbert space with the inner product
	$$\lf(f,g\rt)_{H^q(\mathbf R^N)} = \displaystyle \sum_{|\beta| \leq q}  \int_{\mathbf R^N} \mathscr D^\beta (f) \mathscr D^\beta (g)  dx. $$
	If $q\in \mathbf R$, the inner product and norm of $f,g \in H^q(\mathbf R^N)$ are defined by
	$$\lf(f,g\rt)_{H^q(\mathbf R^N)} = (2 \pi)^N  \int_{\mathbf R^N}  \lf<\xi\rt>^{2q} \widehat{f}(\xi) \overline{\widehat{g}}(\xi) d \xi, \quad \lf<\xi\rt>=(1+|\xi|^2)^{\frac{1}{2}}. $$
	where $\widehat f$ is the Fourier transform of $f$ and $\bar{g}$ is the conjugate term of $g$. Its corresponding norm is 
	$$\lf\|f\rt\|_{H^q(\mathbf R^N)} = (2 \pi)^N  \lf(\displaystyle  \int_{\mathbf R^N}  \lf<\xi\rt>^{2q}  |\widehat{f}(\xi) |^2 d \xi \rt)^{\frac{1}{2}}.$$

	For the above functional space definitions, we refer the interested readers to \cite[Chapter 1, page 11-22]{C003}.
	
	\begin{definition}[Mild solution, see e.g. \cite{OT91,CTG009}]
		If $X$ is a Banach space, we say that $\psi \in C([0,T];X)$ is a mild $X$-value solution of \eqref{Pro} if it satisfies the Duhamel-type integral equation 
		\begin{align} \label{exact-sol}
			\psi(t) = \mathcal Q(t) \psi_0 - \mathrm i \la \int_0^t \mathcal Q(t-\nu) \lf\{|\psi|^{2\eta}(\nu) \psi(\nu) \rt\} d \nu, \quad\mbox{for}~ t \in [0,T],
		\end{align}
		where $\mathcal Q(t)= e^{\mathrm it\Delta}$.
	\end{definition}
	
	Next, we introduce the time behavior of the solution in the three definitions below (see also  \cite{TAX020}).
	\begin{definition}[Continuation] \label{def-cont}
		Given a mild solution $\psi \in L_t^\infty H_x^q (\mathbf B^N)$ of \eqref{Pro} with $\mathbf B^N = (0,T) \times \mathbf R^N$, we say that $\psi_{cont}$ is a  continuation  of $\psi$ if $\psi_{cont} \in L_t^\infty H_x^q(\mathbf B^N_+)$ with  $\mathbf B^N_+ = (0,T_+) \times \mathbf R^N$ is also a mild solution of \eqref{Pro} for some $ T_+ > T$, and $\psi_{cont}( t) = \psi( t)$ whenever $t \in [0,T]$.
	\end{definition}
	\begin{definition}[Maximal time for existence] \label{def-max} Let $\psi (t) $ be the (unique) mild  solution of \eqref {Pro}. We define a maximal time $T_\star $ for existence of the solution $\psi (t) $ to Problem \eqref{Pro} by the following setting:
		\begin{itemize}
			\item [(i)] If $ \psi (t) $ exists for all  $t \in [0 ,\infty)$, then  $T_\star = \infty $.
			\item [(ii)] If there exists a positive constant $T < \infty$ such that the solution $ \psi (t)$ to Problem \eqref{Pro} exists for $t \in [0 ,T)$, but this solution does not exist at $ t = T$, then  $T_\star \equiv T$.
		\end{itemize}
	\end{definition}
	\begin{definition}[Finite-time blow-up] \label{def-blow} Let $\psi(t)$ be a mild solution of \eqref{Pro}. We say $\psi(t)$
		blows up in finite-time if the maximal time $T_\star$ (for existence of Problem \eqref{Pro}) is finite and satisfies
		$
			\lim_{t \uparrow  T_\star^- < \infty} \lf\|\psi(\cdot,t)\rt\|_{H^q(\mathbf R^N)} = \infty.
		$
	\end{definition}
	
	The main tool for our later proofs are the following lemmas. First,  we recall the Gagliardo--Nirenberg interpolation inequality.
	
	\begin{lemma}[Gagliardo--Nirenberg interpolation inequality] \label{GN} Let $s \geq 1, r \leq  \infty$ and $\gamma \in \mathbf N^*$. Suppose also that a real number $\nu$ and $\beta \in \mathbf N^*$ are such that 
		\begin{align} \label{Cond:GN}
			\begin{cases}
				\frac{1}{p} = \frac{\beta}{N} + \lf(\frac{1}{r} - \frac{\gamma}{N}\rt) \nu + \frac{1 - \nu}{s}, \\ \frac{\beta}{\gamma} \leq \nu \leq 1.
			\end{cases}
		\end{align}  
		Then, for every function $f: \mathbf R^N \to \mathbf R$ that lies in $L^s (\mathbf R^N)$ with $\gamma$-th derivative in $L^r (\mathbf R^N)$ also has $\beta$-th derivative in $L^p (\mathbf R^N)$ and there exists a positive constant $C^{\mathrm{GN}}$ which depends on $\beta,\gamma,N, s, r$ and $\nu$, such that
		\begin{align} \label{ineq:GN}
			\lf\|\mathscr D^\beta (f)\rt\|_{L^p (\mathbf R^N)} \leq C^{\mathrm{GN}} \lf\|\mathscr D^\gamma (f)\rt\|_{L^r (\mathbf R^N)}^\nu \lf\|f\rt\|_{L^s (\mathbf R^N)}^{1-\nu}.
		\end{align}
	\end{lemma}
	\noindent \textit{Proof.} See \cite{BM018}, Theorem 1, page 2.  \hfill \rule{1.5mm}{3.5mm} 
	
	Next, we state a version of the Sobolev embedding by the following lemma.
	\begin{lemma}[Sobolev embedding] \label{SE} \label{Sobolev} If $N < 2q  < \infty$, then $H^q(\mathbf R^N) \hookrightarrow C_0(\mathbf R^N)$ and there exists a constant positive $C^{\mathrm{SE}} = C(N,q)$ such that
		$$
		\lf\|f\rt\|_{L^\infty(\mathbf R^N)} \leq C^{\mathrm{SE}} \lf\|f\rt\|_{H^q(\mathbf R^N)}. 
		$$
	\end{lemma}
	\noindent \textit{Proof.} For $N < 2q$, suppose that $f \in S(\mathbf R^N)$ (Schwartz space) we have  
	\begin{align*}
		\lf\|f\rt\|_{L^\infty(\mathbf R^N)} &= \sup_{x \in \mathbf R^N} \lf|\int_{\mathbf R^N} \widehat{f}(\xi) e^{\mathrm i\xi x} d \xi\rt|   \leq \int_{\mathbf R^N}  \lf|\widehat{f}(\xi) \rt|d \xi = \int_{\mathbf R^N}  \frac{1}{\lf<\xi\rt>^{q}}\lf<\xi\rt>^{q} \lf|\widehat{f}(\xi) \rt|d \xi,  \quad \mbox{for}~ \lf<\xi\rt>=(1+|\xi|^2)^{\frac{1}{2}}.
	\end{align*}
	Using the H\"older inequality we conclude that
	\begin{align} \label{Sob}
		\lf\|f\rt\|_{L^\infty(\mathbf R^N)} &\leq \lf(\int_{\mathbf R^N}  \frac{1}{\lf<\xi\rt>^{2q}} d \xi \rt)^{\frac{1}{2}} \lf(\int_{\mathbf R^N}  \lf<\xi\rt>^{2q} \lf|\widehat{f}(\xi) \rt|^2 d \xi\rt)^{\frac{1}{2}}.
	\end{align}
	Since $S(\mathbf R^N) \subseteq H^q(\mathbf R^N)$, this implies that $f \in H^q(\mathbf R^N)$ and $f$ is a Schwartz function, then $f \in C_0(\mathbf R^N)$. Moreover, for $N < 2q$, the first integral in \eqref{Sob} is convergent. Hence, we deduce that
	\begin{align*}
		\lf\|f\rt\|_{L^\infty(\mathbf R^N)} \leq C^{\mathrm{SE}} \lf(\int_{\mathbf R^N} \lf<\xi\rt>^{2q} \lf|\widehat{f}(\xi) \rt|^2 d \xi\rt)^{\frac{1}{2}} \lesssim \lf\|f\rt\|_{H^q(\mathbf R^N)},
	\end{align*}
	where $C = C^{\mathrm{SE}} = \lf(\displaystyle \int_{\mathbf R^N}  \frac{1}{\lf<\xi\rt>^{2q}} d \xi \rt)^{\frac{1}{2}}>0$, for $q > \frac{N}{2}$.
	This concludes the proof. \hfill \rule{1.5mm}{3.5mm}
	
	From now on, we write $\mathcal F_\eta(\psi) = \la |\psi|^{2\eta} \psi$ and assume that $\mathcal F_\eta$ is a smooth function and satisfies the hypotheses:
	\begin{itemize}
		\item For $\eta >0$ and $f,g \in \mathbf R$, there is a Lipschitz constant $C^{\mathrm L}=C(\la,\eta)$ such that
		\begin{align} \label{H1}
			\lf|\mathcal F_\eta(f) - \mathcal F_\eta(g) \rt| \leq C^{\mathrm L} \lf|f-g\rt| \lf(\lf|f\rt|^{2\eta} + \lf|g\rt|^{2\eta} \rt) . \tag{H1}
		\end{align}
		\item There exists an integer number $\sigma$ and $\mathcal F_\eta^{(\sigma)}$ is called the $\sigma-$ th derivative of the function $\mathcal F_\eta$ such that for  $f,g \in \mathbf R$, we have  		
		\begin{align} \label{H2}
			\lf|\mathcal F_\eta^{(\sigma)}(f) - \mathcal F_\eta^{(\sigma)}(g) \rt| \leq C^{\mathrm L} \lf|f-g\rt| \lf(\lf|f\rt|^{2\eta-\sigma} + \lf|g\rt|^{2\eta-\sigma}\rt) . \tag{H2}
		\end{align}
	\end{itemize}
	
	We can now prove our main results.  
	\section{Proof of the main results}

	We shall begin to prove the results in \autoref{existence} - \autoref{Morawetz estimate}, \autoref{prop:space-time-est} and \autoref{decay}. 
	
	\subsection{Proof of Theorem \ref{existence}} We prove the existence and continuation of the solutions to Problem \eqref{Pro} by a fixed point argument. We split the proof into two parts.
	
	\underline{\textit{Part I}: \textit{Local existence}}. Let $N\geq 1$, $\frac{N}{2} < q \in \mathbf N$ and $[1,q] \ni \eta \in \mathbf N$. For positive constants $ K, T $ to be chosen later, we define a closed set with $\mathbf B^N = (0,T) \times \mathbf R^N$ and 
	\begin{align} 
		\mathscr E_K = \lf\{\psi \in L_t^\infty H_x^q (\mathbf B^N) : \psi(\cdot,0) = \psi_0 ~~ \mbox{and} ~~ \lf\|\psi - \psi_0\rt\|_{L_t^\infty H_x^q (\mathbf B^N)} \leq K\rt\},
	\end{align}
	associated with the distance
	$\textrm{dist} (f_1,f_2) = \lf\|f_1-f_2 \rt\|_{L_t^\infty H_x^q(\mathbf B^N)}$ for $f_j \in \mathscr E_K, j=1,2$. We can easily prove that $(\mathscr E_K,\textrm{dist})$ is a complete metric space.
	Let us define the mapping $\mathcal J$ as follows: 
	\begin{align} \label{M}
		\mathcal J \psi (t) = \mathcal Q(t) \psi_0 - \mathrm i \la \int_0^t \mathcal Q(t-\nu) \lf\{|\psi|^{2\eta}(\nu) \psi(\nu) \rt\}  d \nu.
	\end{align} 
	
	Our goal is to prove that $\mathcal J$ is invariant and contraction mapping on $\mathscr E_K$.
	We write \eqref{M} as a fixed-point equation (recall that $\mathcal F_\eta(\psi) = \la |\psi|^{2\eta} \psi$)
	\begin{align}
		&\psi = \mathcal J \psi, \nn\\
		&\mathcal J:  L_t^\infty H_x^q (\mathbf B^N) \to L_t^\infty H_x^q (\mathbf B^N), \nn\\
		&\mathcal J \psi (t) = \mathcal Q(t)  \psi_0 - \mathrm i \int_0^t \mathcal Q(t-\nu) \mathcal F_\eta(\psi) (\nu) d \nu. \label{Phi}
	\end{align}
	We shall write $\mathcal J$ in \eqref{Phi} by
	\begin{align} \label{U-V}
		\mathcal J \psi (t) = \mathcal Q(t)  \psi_0 + \mathcal G (\psi)(t), \quad
		\mathcal G (\psi)(t) = - \mathrm i \int_0^t \mathcal Q(t-\nu) \mathcal F_\eta(\psi) (\nu) d \nu.
	\end{align}
	
	We need to prove that $\mathcal J$ is invariant and a contraction mapping. Indeed, for $\psi_0 \in H^q (\mathbf R^N)$, with $\mathbf N \ni q > \frac{N}{2}$ and $q \geq \eta \geq 1$ and since $\mathcal Q(t)$ is a unitary group, we have
	\begin{align} \label{(e-I)f}
		\lf\|\lf(\mathcal Q(t)  - I\rt) \psi_0 \rt\|^2_{H^q(\mathbf R^N)}  
		&\leq (2\pi)^N \sup_{\xi \in \mathbf R^N} \lf|e^{-\mathrm it|\xi|^2} -1\rt|^2 \int_{\mathbf R^N} \lf<\xi\rt>^{2q} \lf|\widehat{\psi}_0(\xi)\rt|^2 d \xi  \nn\\ 
		&\leq C \lf\|\psi_0\rt\|^2_{H^q(\mathbf R^N)} \lesssim \lf\|\psi_0\rt\|^2_{H^q(\mathbf R^N)}, \quad C=(2\pi)^N>0,
	\end{align}
	where, noting that $\lf|e^{-\mathrm it|\xi|^2} -1\rt| <1$ for all $\xi \in \mathbf R^N$. 
	Let $\psi_1, \psi_2 \in \mathscr E_K$, arguing as in \eqref{(e-I)f}, for every $t \in [0,T]$, we deduce
	\begin{align} \label{Ju-Jv}
		\lf\| \mathcal G (\psi_1)(t) - \mathcal G (\psi_2)(t) \rt\|_{H^q(\mathbf R^N)} \lesssim  \int_0^t \lf\|\mathcal F_\eta(\psi_1) (\nu) - \mathcal F_\eta(\psi_2) (\nu) \rt \|_{H^q(\mathbf R^N)} d \nu  .
	\end{align} 
	For a multi-index $\beta$ with $|\beta|=q \in \mathbf N$,  we have that
	\begin{align*}
		\mathscr D^\beta \lf(\mathcal F_\eta (\psi_1) - \mathcal F_\eta (\psi_2) \rt) = \sum_{\eta \leq q = |\beta|} \lf(\mathcal F_\eta^{(\eta)} (\psi_1)\prod_{j=1}^{\eta} \mathscr D^{\al_j}(\psi_1) - \mathcal F_\eta^{(\eta)} (\psi_2)\prod_{j=1}^{\eta} \mathscr D^{\al_j} (\psi_2) \rt),
	\end{align*}
	where $\eta$ is an integer number such that 
	\begin{align*}
		\begin{cases}
			1 \leq \eta \leq q, \\		
			|\al_j| \geq 1, \quad j = \overline{1,\eta}, \\
			q = |\al_1| + |\al_2| + \cdots + |\al_\eta| = \displaystyle \sum_{j=1}^\eta |\al_j| = |\beta|.
		\end{cases}
	\end{align*}
	Thanks to the triangle inequality,
	\begin{align} \label{[V-V]}
		\lf\|\mathcal F_\eta^{(\eta)} (\psi_1)\prod_{j=1}^{\eta} \mathscr D^{\al_j} (\psi_1) - \mathcal F_\eta^{(\eta)} (\psi_2)\prod_{j=1}^{\eta} \mathscr D^{\al_j}(\psi_2 )\rt\|_{L^2(\mathbf R^N)} \leq (\mbox{Term}_1) + (\mbox{Term}_{2}),
	\end{align}
	where we define 
	$$(\mbox{Term}_{1}) := \lf\|\lf(\mathcal F_\eta^{(\eta)} (\psi_1) - \mathcal F_\eta^{(\eta)} (\psi_2) \rt) \prod_{j=1}^{\eta} \mathscr D^{\al_j} (\psi_1) \rt\|_{L^2(\mathbf R^N)},\quad  (\mbox{Term}_{2})  := \lf\|\mathcal F_\eta^{(\eta)} (\psi_2) \prod_{j=1}^{\eta} \mathscr D^{\al_j}(\psi_1-\psi_2) \rt\|_{L^2(\mathbf R^N)}.$$
	
	We estimate the first in \eqref{[V-V]}. 
	Now let 
	\begin{align} \label{numer}
		p_j = \frac{2q}{|\al_j|}, \quad j = \overline{1,\eta} \quad  \mbox{such that} \quad \frac{2}{p_1} + \frac{2}{p_2} + \cdots + \frac{2}{p_\eta}  = \sum_{j=1}^{\eta} \frac{2}{p_j} = 1,
	\end{align} 
	we continue to use the generalized H\"older inequality  and we have
	\begin{align} \label{est-A}
		(\mbox{Term}_1) \leq  \lf\|\mathcal F_\eta^{(\eta)} (\psi_1) - \mathcal F_\eta^{(\eta)} (\psi_2) \rt\|_{L^\infty(\mathbf R^N)} \prod_{j=1}^{\eta} \lf\|\mathscr D^{\al_j}  (\psi_1) \rt\|_{L^{p_j}  (\mathbf R^N)}.
	\end{align}
	Using hypothesis \eqref{H2} and using the Sobolev embedding (\autoref{SE}) 
	$H^q (\mathbf R^N) \hookrightarrow L^\infty  (\mathbf R^N),$ we have
	\begin{align} \label{u-v}
		\lf\|\mathcal F_\eta^{(\eta)} (\psi_1) - \mathcal F_\eta^{(\eta)} (\psi_2) \rt\|_{L^\infty(\mathbf R^N)} \leq C^{\mathrm L} C^{\mathrm {SE}} \lf(\lf\|\psi_1 \rt\|_{H^q(\mathbf R^N)}^{\eta} + \lf\|\psi_2 \rt\|_{H^q(\mathbf R^N)}^{\eta} \rt) \lf\|\psi_1 - \psi_2 \rt\|_{H^q(\mathbf R^N)}.
	\end{align}
	Next, we estimate the term $\lf\|\mathscr D^{\al_j}  (\psi_1) \rt\|_{L^{p_j} (\mathbf R^N) }$ in \eqref{est-A} by using the Gagliardo-Nirenberg interpolation inequality (\autoref{GN}). Indeed, from the condition \eqref{Cond:GN}, one gets
	\begin{align*}
		\frac{1}{p_j} &= \frac{\al_j}{N} + \lf(\frac{1}{2} - \frac{\gamma_j}{N} \rt) \nu + \frac{1-\nu}{\infty} \nn\\
		&= \frac{\al_j}{N} + \lf(\frac{1}{2} - \frac{\gamma_j}{N} \rt) \nu , \quad \mbox{here, we identify}~ (\al_j, p_j,\gamma_j,\infty) \equiv (\beta, p ,\gamma,s) ~ \mbox{in \autoref{Cond:GN}}.
	\end{align*}
	From this condition, we choose $\nu = \frac{|\al_j|}{q}$, from \eqref{numer} we find that $\gamma_j = q$. Thanks to the Gagliardo-Nirenberg interpolation inequality (\autoref{GN}) we conclude that
	\begin{align*}
		\lf\|\mathscr D^{\al_j} (\psi_1) \rt\|_{L^{p_j} (\mathbf R^N)} \leq C^{\mathrm{GN}} \lf\|\mathscr D^{\gamma_j} (\psi_1) \rt\|_{L^{2}  (\mathbf R^N)}^{\frac{|\al_j|}{q}} \lf\|\psi_1 \rt\|_{L^{\infty} (\mathbf R^N) }^{1-\frac{|\al_j|}{q}} .
	\end{align*}
	From $q=\gamma_j$ yields that $\lf\|\mathscr D^{\gamma_j} (\psi_1) \rt\|_{L^{2}  (\mathbf R^N)} \leq  \lf(\sum_{|\gamma_j| \leq q} \lf\|\mathscr D^{\gamma_j} \psi_1\rt\|_{L^2(\mathbf R^N)}^2\rt)^{\frac{1}{2}} = \lf\|\psi_1 \rt\|_{H^q(\mathbf R^N)}$, we obtain
	\begin{align} \label{D(u)}
		\lf\|\mathscr D^{\al_j} (\psi_1) \rt\|_{L^{p_j} (\mathbf R^N)} \leq C^{\mathrm{GN}} \lf\|\psi_1 \rt\|_{H^q(\mathbf R^N)}^{\frac{|\al_j|}{q}} \lf\|\psi_1 \rt\|_{L^{\infty} (\mathbf R^N) }^{1-\frac{|\al_j|}{q}} .
	\end{align}
	From \eqref{est-A} - \eqref{D(u)}, we deduce that for $C = C^{\mathrm{GN}} C^{\mathrm{SE}} C^{\mathrm{L}} >0$
	\begin{align*} 
		(\mbox{Term}_1) & \leq C^{\mathrm{GN}} C^{\mathrm{SE}} C^{\mathrm{L}} \prod_{j=1}^{\eta} \lf\|\psi_1\rt\|_{H^q(\mathbf R^N)}^{\frac{|\al_j|}{q}} \lf\|\psi_1\rt\|_{L^{\infty} (\mathbf R^N) }^{1-\frac{|\al_j|}{q}} \lf(\lf\|\psi_1 \rt\|_{H^q(\mathbf R^N)}^{\eta} + \lf\|\psi_2 \rt\|_{H^q(\mathbf R^N)}^{\eta} \rt) \lf\|\psi_1 - \psi_2 \rt\|_{H^q(\mathbf R^N)} \nn\\
		& \lesssim \prod_{j=1}^{\eta} \lf\|\psi_1\rt\|_{H^q(\mathbf R^N)}^{\frac{|\al_j|}{q}} \lf\|\psi_1\rt\|_{L^{\infty} (\mathbf R^N) }^{1-\frac{|\al_j|}{q}} \lf(\lf\|\psi_1 \rt\|_{H^q(\mathbf R^N)}^{\eta} + \lf\|\psi_2 \rt\|_{H^q(\mathbf R^N)}^{\eta} \rt) \lf\|\psi_1 - \psi_2 \rt\|_{H^q(\mathbf R^N)}.
	\end{align*}
	Thanks to \autoref{SE}, we have the Sobolev embedding $H^q (\mathbf R^N) \hookrightarrow L^\infty (\mathbf R^N)$ and noting that from \eqref{numer} we derive
	$$\frac{|\al_1|}{q} + \frac{|\al_2|}{q} + \cdots + \frac{|\al_\eta|}{q} = \sum_{j=1}^{\eta} \frac{|\al_j|}{q} = 1,$$
	and for $\eta >1$ $$0 \leq \eta -1 = \eta - \frac{|\al_1|}{q} - \frac{|\al_2|}{q} - \cdots - \frac{|\al_\eta|}{q} = \eta - \frac{|\al_1|+|\al_2|+ \cdots +|\al_\eta|}{q}=  \sum_{j=1}^{\eta} \lf( 1-\frac{|\al_j|}{q} \rt),$$  
	thus we deduce the estimation for $C = C^{\mathrm{GN}} (C^{\mathrm{SE}})^2 C^{\mathrm{L}} >0$
	\begin{align} \label{(u-v)Du}
		(\mbox{Term}_1) \lesssim \lf(\lf\|\psi_1 \rt\|_{H^q(\mathbf R^N)}^{2\eta} + \lf\|\psi_1 \rt\|_{H^q(\mathbf R^N)}^{\eta} \lf\|\psi_2 \rt\|_{H^q(\mathbf R^N)}^{\eta} \rt) \lf\|\psi_1 - \psi_2 \rt\|_{H^q(\mathbf R^N)}. 
	\end{align}
	Next, we estimate $\mbox{Term}_{2}$. Also thanks to the setting for $j = \overline{1,\eta}$ $$p_j = \frac{2q}{|\al_j|} \quad \mbox{and} \quad  \sum_{j=1}^{\eta} \frac{2}{p_j} = 1 ,$$ we continue to use the generalized H\"older inequality  and we obtain 
	\begin{align} \label{est-B}
		&(\mbox{Term}_{2}) \leq \lf\|\mathcal F_\eta^{(\eta)} (\psi_2)  \rt\|_{L^\infty(\mathbf R^N)} \prod_{j=1}^{\eta} \lf\|\mathscr D^{\al_j} (\psi_1-\psi_2) \rt\|_{L^{p_j} (\mathbf R^N)}.
	\end{align}
	Using hypothesis \eqref{H2} and thanks to \autoref{SE}, the Sobolev embedding $H^q(\mathbf R^N) \hookrightarrow L^\infty(\mathbf R^N)$, we deduce that
	\begin{align} \label{est-G}
		\lf\|\mathcal F_\eta^{(\eta)} (\psi_2) \rt\|_{L^\infty(\mathbf R^N)} \leq C^{\mathrm L} C^{\mathrm{SE}} \lf\|\psi_2 \rt\|_{H^q(\mathbf R^N)}^{\eta+1}.
	\end{align}
	From the conditions 
	\begin{align*}
		\frac{1}{p_j}
		= \frac{\al_j}{N} + \lf(\frac{1}{2} - \frac{\gamma_j}{N} \rt) \nu , \quad \nu = \frac{|\al_j|}{q}, \quad \gamma_j =q,
	\end{align*}
	and the Gagliardo-Nirenberg interpolation inequality (\autoref{GN}), we infer that
	\begin{align} \label{D(u-v)}
		\lf\|\mathscr D^{\al_j} 		(\psi_1-\psi_2)\rt\|_{L^{p_j} (\mathbf R^N)} \leq C^{\mathrm{GN}} \lf\|\psi_1-\psi_2\rt\|_{H^q(\mathbf R^N)}^{\frac{|\al_j|}{q}} \lf\|\psi_1-\psi_2\rt\|_{L^{\infty} (\mathbf R^N) }^{1-\frac{|\al_j|}{q}} .
	\end{align}
	From \eqref{est-B} - \eqref{D(u-v)}, using the Sobolev embedding $H^q (\mathbf R^N) \hookrightarrow L^\infty (\mathbf R^N)$ and noting that 
	$$\sum_{j=1}^{\eta} \frac{|\al_j|}{q} = 1, \quad \mbox{and} ~~0 \leq \eta-1 = \sum_{j=1}^{\eta} \lf( 1-\frac{|\al_j|}{q} \rt),$$ 
	we have that, for $C = C^{\mathrm{GN}}(C^{\mathrm{SE}})^2 C^{\mathrm{L}} > 0$,
	\begin{align} \label{VuD(u-v)}
		(\mbox{Term}_{2}) &\leq C^{\mathrm{GN}} C^{\mathrm{SE}} C^{\mathrm{L}} \lf\|\psi_2 \rt\|_{H^q(\mathbf R^N)}^{\eta+1} \prod_{j=1}^{\eta} \lf\|\psi_1-\psi_2\rt\|_{H^q(\mathbf R^N)}^{\frac{|\al_j|}{q}} \lf\|\psi_1-\psi_2\rt\|_{L^{\infty} (\mathbf R^N) }^{1-\frac{|\al_j|}{q}} \nn\\
		&\leq C^{\mathrm{GN}}(C^{\mathrm{SE}})^2 C^{\mathrm{L}} \lf(\lf\|\psi_2 \rt\|_{H^q(\mathbf R^N)}^{2\eta} + \lf\|\psi_2 \rt\|_{H^q(\mathbf R^N)}^{\eta+1} \lf\|\psi_1 \rt\|_{H^q(\mathbf R^N)}^{\eta-1} \rt) \lf\|\psi_1-\psi_2\rt\|_{H^q(\mathbf R^N)} 
		\nn\\
		&\lesssim  \lf(\lf\|\psi_2 \rt\|_{H^q(\mathbf R^N)}^{2\eta} + \lf\|\psi_2 \rt\|_{H^q(\mathbf R^N)}^{\eta+1} \lf\|\psi_1 \rt\|_{H^q(\mathbf R^N)}^{\eta-1} \rt) \lf\|\psi_1-\psi_2\rt\|_{H^q(\mathbf R^N)}. 
	\end{align}
	For $\eta\geq 1$, combining \eqref{[V-V]}, \eqref{(u-v)Du} and \eqref{VuD(u-v)}, we have
	\begin{align*}  
		\lf\|\mathcal F_\eta(\psi_1) - \mathcal F_\eta(\psi_2)  \rt\|_{H^q(\mathbf R^N)}  & = \lf(\sum_{|\beta| \leq q} \lf\|\mathscr D^\beta \lf(\mathcal F_\eta(\psi_1) - \mathcal F_\eta(\psi_2) \rt) \rt\|_{L^2(\mathbf R^N)}^2\rt)^{\frac{1}{2}} \nn \\
		& \lesssim \lf(\sum_{|\beta| \leq q} \sum_{\eta \leq q} \lf\|\mathcal F_\eta^{(\eta)} (\psi_1)\prod_{j=1}^{\eta} \mathscr D^{\al_j} (\psi_1) - \mathcal F_\eta^{(\eta)} (\psi_2)\prod_{j=1}^{\eta} \mathscr D^{\al_j} (\psi_2) \rt\|_{L^2(\mathbf R^N)}^2\rt)^{\frac{1}{2}} \nn\\ 
		&\lesssim L(K, \lf\|\psi_0\rt\|_{H^q(\mathbf R^N)}) \lf\|\psi_1-\psi_2\rt\|_{H^q(\mathbf R^N)},  
	\end{align*}
	where $L$ is dependent on $K$ and $\lf\|\psi_0\rt\|_{H^q(\mathbf R^N)}$.
	From this and \eqref{Ju-Jv}, we derive 
	\begin{align} \label{Mu-Mv}
		\lf\|\mathcal J\psi_1(t) - \mathcal J\psi_2(t) \rt\|_{H^q(\mathbf R^N)}  
		&\lesssim T L(K, \lf\|\psi_0\rt\|_{H^q(\mathbf R^N)}) \lf\|\psi_1-\psi_2 \rt\|_{L_t^\infty H_x^q(\mathbf B^N)}. 
	\end{align}
	From \eqref{Ju-Jv} and \eqref{Mu-Mv}, if we take $\psi \equiv \psi_1 \in \mathscr E_K$ and $\psi_2 = 0$, it follows 
	\begin{align} \label{J_psi}
		\lf\|\mathcal G(\psi)  \rt\|_{L_t^\infty H_x^q (\mathbf B^N) } &\lesssim T L (K, \lf\|\psi_0\rt\|_{H^q(\mathbf R^N)}) \lf\|\psi \rt\|_{L_t^\infty H_x^q (\mathbf B^N) } \nn\\
		&\lesssim T L (K, \lf\|\psi_0\rt\|_{H^q(\mathbf R^N)}) (K + \lf\|\psi_0 \rt\|_{H^q(\mathbf R^N)} ) = T L_1 (K, \lf\|\psi_0\rt\|_{H^q(\mathbf R^N)}) ,
	\end{align}
	where $L_1 (K, \lf\|\psi_0\rt\|_{H^q(\mathbf R^N)}) = L (K, \lf\|\psi_0\rt\|_{H^q(\mathbf R^N)})  (K + \lf\|\psi_0 \rt\|_{H^q(\mathbf R^N)} )$. Hence, if we choose $K$ such that $K = 2C\lf\|\psi_0\rt\|_{H^q(\mathbf R^N)}$ and 
	$2C T L_1 (K, \lf\|\psi_0\rt\|_{H^q(\mathbf R^N)})   < K$, together with \eqref{U-V} and \eqref{(e-I)f}, to obtain
	\begin{align*}
		\lf\|\mathcal J\psi \rt\|_{L_t^\infty H_x^q(\mathbf B^N)} &\leq \lf\|\lf(\mathcal Q  - I\rt) \psi_0 \rt\|_{L_t^\infty H_x^q(\mathbf B^N) } + \lf\|\mathcal G(\psi)  \rt\|_{L_t^\infty H_x^q(\mathbf B^N) } \nn\\
		& \lesssim \lf\|\psi_0\rt\|_{H^q(\mathbf R^N)} + T L_1 (K, \lf\|\psi_0\rt\|_{H^q(\mathbf R^N)}) \leq K,
	\end{align*} 
	then $\mathcal J$ is an invariant mapping. We next prove that $\mathcal J$ is a contraction. From \eqref{Mu-Mv}, the right-hand side is independent of $t$, we deduce that 
	$\textrm{dist} \lf(\mathcal J \psi_1, \mathcal J \psi_2 \rt) \lesssim \rho  \textrm{dist}  \lf(\psi_1,\psi_2 \rt)$
	for $\rho := CT L(K, \lf\|\psi_0\rt\|_{H^q(\mathbf R^N)}) >0$.
	Choosing $T$ and $K$ small enough such that $\rho < 1$, then the mapping $\mathcal J$ is a contraction on $\mathscr E_K$. Hence, by using the contraction mapping principle, we conclude that
	the mapping $\mathcal J$ admits a unique fixed point $\psi$ in $\mathscr E_K$.
	
	\underline{\textit{Part  II:}  Continuation of solutions}. Let $\psi: [0,T] \to H^q(\mathbf R^N),~ (\mathbf N \ni q > \frac{N}{2}$, $[1,q] \ni \eta \in \mathbf N$) be a unique mild solution of Problem \eqref{Pro}, and $T$ be the positive time for existence (defined as in Part I). Fix $K >0,$ and for a larger time $T_+> T$ ($T_+$ depends on $K$),  we shall prove that $\psi_{cont}: [0,T_+] \to H^q(\mathbf R^N)$ is also the unique mild solution of Problem \eqref{Pro}. This should be well dealt with $4\max\{\mathcal R_j\} \leq K$, $j=\overline{1,3}$ where the definitions for $C>0$ are
	\begin{subequations}
		\begin{align}
			&\mathcal R_1 := C \lf\|\psi_0\rt\|_{ H^q (\mathbf R^N)}; \label{I1} \\
			&\mathcal R_2 := CT_+ L_1 (K, \lf\|\psi(\cdot,T)\rt\|_{H^q(\mathbf R^N)}); \label{I 2} \\
			&\mathcal R_3 := C L (K,\lf\|\psi(\cdot,T)\rt\|_{H^q(\mathbf R^N)}) T_+. \label{I 4} 
		\end{align}
	\end{subequations}
	Here, $L$, $L_1$ are the Lipschitz constants as in \eqref{Mu-Mv} and \eqref{J_psi}.  
	For $T_+ \geq T >0$ (noting that one can choose $T_+ >T$ and close enough to $T$) and $K>0$, let us define for $\mathbf B^N_+ = (0,T_+) \times \mathbf R^N$
	\begin{align} \label{wide-U}
		\widetilde{\mathscr E}_{K} =\lf\{ \psi_{cont} \in L_t^\infty H_x^q (\mathbf B^N_+)  \bigg|\hspace*{-0.4cm}\begin{array}{lll}
			&\psi_{cont} (t) =\psi(t),&\forall t\in [0, T], \\
			&\lf\|\psi_{cont} - \psi(T) \rt\|_{L_t^\infty H_x^q (\mathbf B^N_+)} \leq K, & \forall  t \in [T,T_+], \end{array}
		\rt\}
	\end{align}
	equipped with the distance 
	$\textrm{dist}^+ (f_1,f_2) = \lf\|f_1-f_2 \rt\|_{L_t^\infty H_x^q (\mathbf B^N_+)}$ for $f_i \in \widetilde{\mathscr E}_{K}, i=1,2,$
	and we note that $(\widetilde{\mathscr E}_{K},\textrm{dist}^+)$ is a complete metric space. Next, we need to clarify the following two claims.
	
	$\bullet$ \underline{\textit{Claim 1}}: \textit{We show that $\mathcal J$  defined as in \eqref{M}
		is an operator on $\widetilde {\mathscr E}_{K}$}. Indeed, let $\psi_{cont} \in \widetilde {\mathscr E}_{K}$, we consider the following two cases. \\
	$*$ For $t\in [0,T]$, then by virtue of  Part I, we know that Problem \eqref{Pro} admits a unique solution which satisfies $\mathcal J\psi_{cont} = \psi_{cont}$ and $ \psi_{cont}(\cdot,t)=\psi(\cdot,t)$, thus 		
	$\lf\|\mathcal J  \psi_{cont}  -  \psi \rt\|_{L_t^\infty H_x^q (\mathbf B^N_+)} = 0$ in $\widetilde {\mathscr E}_{K}$ for all $t \in [0,T]$.\\
	$*$ For $t\in [T,T_+]$, from \eqref{exact-sol} and \eqref{M}, we have 
	\begin{align} \label{3N}
		\mathcal J \psi_{cont} (t)-\psi(T)
		& = \lf(\mathcal Q(t) - \mathcal Q(T) \rt)\psi_0  - \mathrm i\int_T^t \mathcal Q(t-\nu) \mathcal F_\eta(\psi_{cont})(\nu) d \nu  \nn\\
		&\quad - \mathrm i\int_0^T \lf(\mathcal Q(t - \nu) - \mathcal Q(T - \nu) \rt)  \mathcal F_\eta(\psi_{cont}) (\nu) d \nu  = \sum_{j=3}^5 (\mbox{Term}_{j}),
	\end{align}
	where we define 
	$$(\mbox{Term}_{3}) := \lf(\mathcal Q(t) - \mathcal Q(T) \rt)\psi_0,$$
	$$(\mbox{Term}_{4}) := - \mathrm i\int_T^t \mathcal Q(t-\nu) \mathcal F_\eta(\psi_{cont})(\nu) d \nu,$$
	$$(\mbox{Term}_{5}) :=- \mathrm i\int_0^T \lf(\mathcal Q(t - \nu) - \mathcal Q(T - \nu) \rt)  \mathcal F_\eta(\psi_{cont}) (\nu) d \nu.$$
	To estimate the first term in the RHS of \eqref{3N}, for all $t \in [T,T_+]$, we have
	\begin{align*} 
		\lf\|(\mbox{Term}_{3})\rt\|_{ H^q (\mathbf R^N)}
		&\lesssim \lf\|\psi_0\rt\|_{H^q(\mathbf R^N)}.
	\end{align*}
	By \eqref{I1}, the following estimate holds
	\begin{align} \label{J1}
		\lf\|(\mbox{Term}_{3})\rt\|_{L_t^\infty H_x^q (\mathbf B^N_+)} \leq \mathcal R_1 \leq  \frac{K}{4}, \quad \forall t \in [0,T_+].
	\end{align}
	From \eqref{wide-U}, we deduce for all $t \in [T,T_+]$  that 
	$$\lf\|\psi_{cont}(\cdot,t)\rt\|_{H^q(\mathbf R^N)} \leq K +  \lf\|\psi(\cdot,T)\rt\|_{H^q(\mathbf R^N)},$$
	and similar to \eqref{J_psi} for all $t \in [T,T_+]$, we have the following estimate
	\begin{align*} 
		\lf\|(\mbox{Term}_{4})\rt\|_{H^q(\mathbf R^N)} \lesssim T_+ L_1 (K, \lf\|\psi (\cdot,T)\rt\|_{H^q(\mathbf R^N)}).
	\end{align*}
	Using \eqref{I 2}, we infer that
	\begin{align}  \label{J2}
		\lf\|(\mbox{Term}_{4})\rt\|_{L_t^\infty H_x^q (\mathbf B^N_+)} \leq \mathcal R_2 \leq \frac{K}{4}.
	\end{align}
	We continue the procedure with estimates for the third term in the RHS of \eqref{3N}, similar to \eqref{J_psi} for all $t \in [T,T_+]$, we obtain that
	\begin{align*} 
		\lf\|(\mbox{Term}_{5})\rt\|_{H^q(\mathbf R^N)} \lesssim T_+ L_1 (K, \lf\|\psi(\cdot,T)\rt\|_{H^q(\mathbf R^N)}).
	\end{align*}
	From \eqref{I 2}, we have that
	\begin{align} \label{J3}
		\lf\|(\mbox{Term}_{5})\rt\|_{L_t^\infty H_x^q (\mathbf B^N_+)} \leq \mathcal R_2 \leq \frac{K}{4}.
	\end{align}
	It follows from \eqref{J1}-\eqref{J3} for all $t\in [0,T_+]$ that
	\begin{align*}
		\lf\|\mathcal J \psi_{cont} - \psi (T) \rt\|_{L_t^\infty H_x^q (\mathbf B^N_+)} \leq \frac{K}{4} + \frac{K}{4} + \frac{K}{4} = \frac{3K}{4} < K.
	\end{align*}%
	Hence, we have shown that $\mathcal J$ is a map from $\widetilde {\mathscr E}_{K}$ to $ \widetilde {\mathscr E}_{K}$.
	
	$\bullet$ \underline{\textit{Claim 2}}: \textit{We show that $\mathcal J$ is a contraction on $\widetilde{\mathscr E}_{K}$}. Let $
	\psi_1,\psi_2 \in \widetilde{\mathscr E}_{K}$, and we have that, for $t \in [0, T_+]$,
	\begin{align*}
		\mathcal J \psi_1 (t) - \mathcal J \psi_2  (t) = \mathrm i \int_0^t  \mathcal Q(t-\nu) \lf(\mathcal F_\eta(\psi_2)(\nu)  - \mathcal F_\eta(\psi_1) (\nu) \rt) d \nu.
	\end{align*}
	For all $t \in [0,T]$,  we conclude that $\mathcal J\psi_1(t)-\mathcal J\psi_2(t) ~~ \mbox{vanishes on} ~~ [0,T]$. For $t \in [T,T_+]$, proceeding as in \eqref{Mu-Mv}, we obtain
	\begin{align*}
		\lf\|\mathcal J\psi_1(t) - \mathcal J\psi_2(t) \rt\|_{H^q(\mathbf R^N)} \lesssim L(K, \lf\|\psi (\cdot,T) \rt\|_{H^q(\mathbf R^N)}) T_+ \lf\|\psi_1-\psi_2 \rt\|_{L_t^\infty H_x^q (\mathbf B^N_+)},
	\end{align*}
	where $L$ is the Lipschitz constant as in \eqref{Mu-Mv}.  
	Hence, from \eqref{I 4}, we deduce that
	\begin{align*}
		\lf\|\mathcal J \psi_1 - \mathcal J \psi_2 \rt\|_{L_t^\infty H_x^q (\mathbf B^N_+)} &\lesssim  L(K, \lf\|\psi(\cdot,T) \rt\|_{H^q(\mathbf R^N)}) T_+ \lf\|\psi_1 - \psi_2 \rt\|_{L_t^\infty H_x^q (\mathbf B^N_+)} \nn\\
		&\leq \mathcal R_3  \leq \frac{K}{4} \lf\|\psi_1-\psi_2\rt\|_{L_t^\infty H_x^q (\mathbf B^N_+)}.
	\end{align*}
	Thus, for any $T_+>0$ and without loss of generality, we may choose the radius $K$ such that $0\leq K <4$, and one obtains that
	\begin{align} \label{EE1}
		\textrm{dist}^+ (\mathcal J \psi_1 , \mathcal J \psi_2 ) \leq \frac{K}{4} \textrm{dist}^+ (\psi_1,\psi_2).
	\end{align}
	This implies that $\mathcal J$ is a $\frac{K}{4}$-contraction. By applying the Banach contraction principle, it follows that $\mathcal J$ has a unique fixed point $\psi_{cont}$ in the ball $\widetilde {\mathscr E}_{K}$, and we can say that $\psi_{cont}$ is a continuation in the larger time of $\psi$. \hfill \rule{1.5mm}{3.5mm}

	\subsection{Proof of Theorem \ref{blow-up}} 
	We shall begin to prove two results in \autoref{blow-up}.
	
	\underline{\textit{Part I}: \textit{Finite time blow-up or global existence}}. For $N \geq 1$, let $\psi_0 \in H^q(\mathbf R^N)$ with $\mathbf N \ni q > \frac{N}{2}$, $\mathbf N \ni \eta \in [1,q]$ and define 
	\begin{equation*}
		T_\star:=\sup\lf\{T>0:\mbox{there exits a solution on}\; [0,T]\rt\}.
	\end{equation*}
	Suppose that $T_\star < \infty$,  and $\|\psi(\cdot,t) \|_{H^q(\mathbf R^N)} \leq K_0,$ for some $K_0 >0$, and for all $t \in (0,T_\star)$. Now assume that there exists a sequence $\{t_{n}\}_{n\in \mathbf N^*}\subset (0,T_\star)$ such that $t_n \uparrow T_\star$. We will show that $\{\psi( t_{n})\}_{n\in \mathbf N^*}$ is a Cauchy sequence in $H^q(\mathbf R^N)$. In fact, given $\ep >0$, fix $\mathcal N_\ep = \mathcal N(\ep) \in \mathbf N^*$ such that for all $n,j>\mathcal N_\ep$: $t_n ,t_{j}\in (0,T_\star)$, and $t_n < t_j$, we have
	\begin{align*}
		\psi( t_{j})-\psi( t_n) &= \lf(\mathcal Q (t_j-t_n) - I \rt)\mathcal Q(t_n) \psi_0  - \mathrm i \int_{t_n}^{t_j} \mathcal Q(t_j - \nu ) \mathcal F_\eta(\psi ) (\nu) d\nu  \\
		&\quad -  \mathrm i \lf(\mathcal Q(t_j - t_n) - I\rt) \int_0^{t_n} \mathcal Q(t_n - \nu ) \mathcal F_\eta(\psi) (\nu) d\nu \nn\\
		&= \lf(\mathcal Q(t_j-t_n) - I \rt) \psi( t_n) - \mathrm i \int_{t_n}^{t_j} \mathcal Q(t_j - \nu ) \mathcal F_\eta(\psi) (\nu)  d \nu.
	\end{align*}
	Thus, we have
	\begin{align*}
		\lf\|\psi(\cdot,t_{j})-\psi(\cdot,t_n)\rt\|_{H^q(\mathbf R^N)} &\leq \underbrace{\lf\|\lf(\mathcal Q(t_j-t_n) - I \rt) \psi( t_n) \rt\|_{H^q(\mathbf R^N)} }_{=:\text{(Term}_{6})} + \underbrace{\lf\| \mathrm i \int_{t_n}^{t_j} \mathcal Q(t_j - \nu ) \mathcal F_\eta(\psi) (\nu) d \nu \rt\|_{H^q(\mathbf R^N)} }_{=:(\text{Term}_{7})} .
	\end{align*}
	First, we estimate the first term in the inequality above, and we have that
	\begin{align*}
		(\text{Term}_{6}) \leq \lf\|\mathcal Q(t_j-t_n) - I  \rt\|_{\mathbb L(H^q(\mathbf R^N) \to H^q(\mathbf R^N))} \lf\|\psi(\cdot,t_n) \rt\|_{H^q(\mathbf R^N)}. \end{align*}
	By an argument analogous to that used in \eqref{J_psi}, we obtain 
	\begin{align*}
		(\text{Term}_{7}) \lesssim  L (K_0, \lf\|\psi_0\rt\|_{H^q(\mathbf R^N)}) K_0 |t_j - t_n|.
	\end{align*}
	Since $\{t_n\}_{n\in\mathbf N^*}$ is a convergent sequence, then we can take $\mathcal N_\ep \in \mathbf N^*$ with
	$j \geq n \geq \mathcal N_\ep$ such that $|t_{j} - t_n|$ is as small as we want. Since the semigroup $\lf\{\mathcal Q(t) \rt\}_{t \geq 0}$ is strongly continuous in $H^q(\mathbf R^N)$, given $\ep>0$, we have
	$$
	K_0 \lf\|\lf(\mathcal Q(t_j-t_n) - I \rt)\rt\|_{\mathbb L(H^q(\mathbf R^N) \to H^q(\mathbf R^N))}  < \frac{\ep}{2},  
	$$
	and
	$$
	L (K_0, \lf\|\psi_0\rt\|_{H^q(\mathbf R^N)}) K_0 |t_j - t_n| < \frac{\ep}{2}.
	$$
	From this, given $\ep >0$, there exists $\mathcal N_\ep \in \mathbf N^*$ such that
	\begin{align*}
		\lf\| \psi(\cdot,t_{j})-\psi(\cdot,t_n)\rt\|_{H^q(\mathbf R^N)} < \ep, \quad \mbox{for}~~j,n \geq \mathcal N_\ep. 
	\end{align*}
	It follows that $\{\psi(\cdot,t_n)\}_{n\in \mathbf N^*} \subset H^q(\mathbf R^N)$
	is a Cauchy sequence.  Hence, for $\{t_n\}_{n\in \mathbf N^*}$, an arbitrary sequence of $(0,T_\star)$, the existence of the following 
	limit has been proved $\lim_{t \uparrow T_\star^- < \infty} \lf\|\psi(\cdot,t)\rt\|_{H^q(\mathbf R^N)} < \infty$.  From the results of \autoref{existence}, we have argued that the existence of the solution can be expanded to some larger time interval  (i.e., $\psi$ may be continued beyond the maximal time $T_\star$). This contradicts the definition of $T_\star$. Therefore, either $ T_\star = \infty$ or if $T_\star<\infty$ then $\lim_{t \uparrow T_\star^-} \lf\|\psi(\cdot,t)\rt\|_{H^q(\mathbf R^N)} = \infty.$
	
	\underline{\textit{Part II}: \textit{Continuous dependence}}. For $\frac{N}{2} < q \in \mathbf N$, $[1,q] \ni \eta \in \mathbf N$, let $\psi_0 \in H^q(\mathbf R^N) $ and consider $\psi_{0j} \in H^q(\mathbf R^N)$ such that
	$\psi_{0j} \to \psi_0 ~~ \mbox{as} ~~ j \uparrow \infty.$
	For $j$ sufficiently large we have that
	\begin{align} \label{con1}
		\lf\|\psi(\cdot,t)-\psi_j(\cdot,t)\rt\|_{H^q(\mathbf R^N)} \leq  \lf\|\mathcal Q(t) \lf(\psi_0 - \psi_{0j}\rt)\rt\|_{H^q(\mathbf R^N)}  + \int_0^t \lf\|\mathcal Q(t-\nu)  \lf( \mathcal F_\eta(\psi)(\nu) - \mathcal F_\eta(\psi_j)(\nu)\rt) \rt\|_{H^q(\mathbf R^N)} d\nu.  
	\end{align}
	Similar to the computations above, we obtain
	\begin{align} \label{con2}
		\lf\|\mathcal Q(t) \lf( \psi_0 - 	\psi_{0j}\rt)\rt\|_{H^q(\mathbf R^N)}
		\lesssim \lf\|\psi_0-\psi_{0j} \rt\|_{H^q(\mathbf R^N)},
	\end{align}
	and we claim that
	\begin{align} \label{con3}
	\int_0^t \lf\|\mathcal Q(t-\nu) \lf(\mathcal F_\eta(\psi)(\nu) - \mathcal F_\eta(\psi_j)(\nu)\rt) \rt\|_{H^q(\mathbf R^N)} d\nu \lesssim  L(K_0, \lf\|\psi_0 \rt\|_{H^q(\mathbf R^N)}) \int_0^t \lf\|\psi( \nu) - \psi_j( \nu)\rt\|_{H^q(\mathbf R^N)} d\nu .  
	\end{align}
	Combining \eqref{con1}-\eqref{con3}, we deduce that 
	\begin{align*} 
		\lf\|\psi( t)-\psi_j( t)\rt\|_{H^q(\mathbf R^N)} \lesssim \lf\|\psi_0-\psi_{0j} \rt\|_{H^q(\mathbf R^N)} + L(K_0, \lf\|\psi_0 \rt\|_{H^q(\mathbf R^N)}) \int_0^t \lf\|\psi( \nu) - \psi_j( \nu) \rt\|_{H^q(\mathbf R^N)} d\nu.
	\end{align*}
	By using the Gr\"onwall inequality, we see that
	\begin{align} \label{con6}
		\lf\|\psi(\cdot, t)-\psi_j(\cdot, t)\rt\|_{H^q(\mathbf R^N)} \lesssim  \lf\|\psi_0-\psi_{0j} \rt\|_{H^q(\mathbf R^N)}  \exp\lf(CL (K_0, \lf\|\psi_0 \rt\|_{H^q(\mathbf R^N)}) t\rt), \quad  t \in [0,T_\star].
	\end{align}
	As $j \uparrow \infty$ we have that $\psi_{0j} \to \psi_0$ then $\psi_j \to \psi$ in $H^q(\mathbf R^N)$. Iterating this property to cover any compact subset of $[0,T_\star]$, i.e., starting from $\psi_0$ and $\psi( t_1)$ over the intervals $[0,t_1]$ and $[t_1,t_2],$ respectively, the proof of \autoref{blow-up} is finished. \hfill \rule{1.5mm}{3.5mm}

\subsection{Proof of Theorem \ref{stability}} As shown above, the solutions of Problem \eqref{Pro} exist and satisfy the conservation laws \eqref{ma}, \eqref{mass} and \eqref{Mom}. We consider the $(2-2)$ Morawetz action
\begin{align} \label{V}
	\mathcal V_{\delta}(t) = \iint_{\mathbf R^N \times \mathbf R^N} |\psi(y,t) |^2  \vec{\slashed{\nabla}}_x \delta(x-y) \cdot \overrightarrow{Q^\nabla} (\psi)(x,t) dx dy.
\end{align}
From now on, when no confusion can arise, we denote $\displaystyle \iint f dx = \iint_{\mathbf R^N \times \mathbf R^N} f dx$. For the mass density $\varPsi = |\psi|^2$, the derivative of both sides with respect to the variable $t$ of the quantity \eqref{V} is
\begin{align} \label{est:K}
	\frac{d}{dt} \mathcal V_\delta(t) =&\underbrace{\iint_{\mathbf R^N \times \mathbf R^N} \frac{\partial \varPsi}{\partial t} (y,t) \vec{\slashed{\nabla}}_x \delta(x-y) \cdot \overrightarrow{Q^\nabla} (\psi)(x,t) dx  dy}_{=:(\textrm{Integral}_1)(t)} + \underbrace{\iint_{\mathbf R^N \times \mathbf R^N} \varPsi(y,t)   \vec{\slashed{\nabla}}_x \delta(x-y) \cdot \frac{\partial \overrightarrow{Q^\nabla} (\psi)}{\partial t} (x,t) dx dy }_{=:(\textrm{Integral}_2)(t)}.
\end{align}
First, from \eqref{mass} one has $\frac{\partial \varPsi}{\partial t} (y,t) = - 2 \mbox{div} \lf(\overrightarrow{Q^\nabla}(\psi)(y,t) \rt)$, thus we have
\begin{align} \label{est:Int1}
	(\textrm{Integral}_1)(t) = -2 \iint_{\mathbf R^N \times \mathbf R^N}  \lf(  \vec{\slashed{\nabla}}_{y} \cdot \overrightarrow{Q^\nabla} (\psi)(y,t) \rt) \lf(\vec{\slashed{\nabla}}_x \delta(x-y) \cdot \overrightarrow{Q^\nabla} (\psi)(x,t) \rt) dx  dy.
\end{align}
Second, for $\vec{\slashed{\nabla}}_x \delta = (\slashed{\nabla}_{x,1} \delta,...,\slashed{\nabla}_{x,N} \delta)$ and $\overrightarrow{Q^\nabla} (\psi) = (Q^\nabla_1 (\psi),...,Q^\nabla_N (\psi))$, we have
\begin{align*}
	(\textrm{Integral}_2)(t) 
	&= \iint_{\mathbf R^N \times \mathbf R^N} \varPsi(y,t)  \slashed{\nabla}_{x,1} \delta(x-y) \frac{\partial Q^\nabla_1 (\psi)}{\partial t} (x,t)  dx dy \\
	&\quad + \iint_{\mathbf R^N \times \mathbf R^N} \varPsi(y,t) \slashed{\nabla}_{x,2} \delta(x-y)  \frac{\partial Q^\nabla_2 (\psi)}{\partial t} (x,t)  dx dy 	\\ 
	&\quad \stackrel{\cdots}{+} \iint_{\mathbf R^N \times \mathbf R^N} \varPsi(y,t) \slashed{\nabla}_{x,N-1} \delta(x-y)  \frac{\partial Q^\nabla_{N-1} (\psi)}{\partial t}  (x,t) dx dy \\
	&\quad + \iint_{\mathbf R^N \times \mathbf R^N} \varPsi(y,t) \slashed{\nabla}_{x,N} \delta(x-y) \frac{\partial Q^\nabla_N (\psi)}{\partial t} (x,t) dx dy.
\end{align*}
From the local momentum conservation \eqref{Mom}, we have for $m = 1,...,N$
\begin{align*} 
	\frac{\partial Q^\nabla_m(\psi)}{\partial t} (x,t) &= -  \slashed{\nabla}_{x,1} \lf(\mathbf 1_{m ,1} L^\nabla(\psi) + S^\nabla_{m, 1}(\psi) \rt) (x,t)  - \slashed{\nabla}_{x,2} \lf(\mathbf 1_{m, 2} L^\nabla(\psi)  + S^\nabla_{m, 2}(\psi)  \rt) (x,t) \nn\\
	&\quad \stackrel{\cdots}{-} \slashed{\nabla}_{x,N-1} \lf(\mathbf 1_{m, N-1} L^\nabla(\psi)  + S^\nabla_{m ,N-1} (\psi) \rt) (x,t) - \slashed{\nabla}_{x,N} \lf(\mathbf 1_{m, N} L^\nabla(\psi)  + S^\nabla_{m, N} (\psi) \rt)(x,t) ,
\end{align*}
and from the properties of Kronecker delta function \eqref{Kronecker}, this leads to the following conclusion
\begin{align*} 
	(\textrm{Integral}_2)(t)
	&=- [\iint_{\mathbf R^N \times \mathbf R^N}  \cdots \iint_{\mathbf R^N \times \mathbf R^N} ] [\varPsi(y,t) \slashed{\nabla}_{x,n} \delta(x-y)    \slashed{\nabla}_{x,n} L^\nabla(\psi) (x,t)]_{n=1}^N dxdy \\
	&\quad - [\iint_{\mathbf R^N \times \mathbf R^N}  \cdots \iint_{\mathbf R^N \times \mathbf R^N} ] [\varPsi(y,t) \slashed{\nabla}_{x,1} \delta(x-y)    \slashed{\nabla}_{x,n}  S^\nabla_{1,n}(\psi) (x,t)]_{n=1}^N dxdy \\
	&\quad - [\iint_{\mathbf R^N \times \mathbf R^N}  \cdots \iint_{\mathbf R^N \times \mathbf R^N} ] [\varPsi(y,t) \slashed{\nabla}_{x,2} \delta(x-y)    \slashed{\nabla}_{x,n}  S^\nabla_{2,n}(\psi) (x,t)]_{n=1}^N dxdy \\
	&\quad \vdots \stackrel{\cdots}{-}	\\ 
	&\quad - [\iint_{\mathbf R^N \times \mathbf R^N}  \cdots \iint_{\mathbf R^N \times \mathbf R^N} ] [\varPsi(y,t) \slashed{\nabla}_{x,N-1} \delta(x-y)    \slashed{\nabla}_{x,n}  S^\nabla_{N-1,n}(\psi) (x,t)]_{n=1}^N dxdy \\
	&\quad - [\iint_{\mathbf R^N \times \mathbf R^N}  \cdots \iint_{\mathbf R^N \times \mathbf R^N} ] [\varPsi(y,t) \slashed{\nabla}_{x,N} \delta(x-y)    \slashed{\nabla}_{x,n}  S^\nabla_{N,n}(\psi) (x,t)]_{n=1}^N dxdy .
\end{align*}
Next, integrating by parts, one has
\begin{align*}
	(\textrm{Integral}_2) (t)
	&= [\iint_{\mathbf R^N \times \mathbf R^N}  \cdots \iint_{\mathbf R^N \times \mathbf R^N} ] [\varPsi(y,t) (\slashed{\nabla}_{x,n}^2 \delta(x-y) )   L^\nabla(\psi) (x,t)]_{n=1}^N  dxdy \\
	&\quad + [\iint_{\mathbf R^N \times \mathbf R^N}  \cdots \iint_{\mathbf R^N \times \mathbf R^N} ] [\varPsi(y,t) \lf( \slashed{\nabla}_{x,1} \slashed{\nabla}_{x,n} \delta(x-y) \rt)   S^\nabla_{1,n} (\psi) (x,t)]_{n=1}^N dxdy \\
	&\quad + [\iint_{\mathbf R^N \times \mathbf R^N}  \cdots \iint_{\mathbf R^N \times \mathbf R^N} ] [\varPsi(y,t) \lf( \slashed{\nabla}_{x,2} \slashed{\nabla}_{x,n} \delta(x-y) \rt)   S^\nabla_{2,n} (\psi) (x,t)]_{n=1}^N dxdy \\
	&\quad \vdots \stackrel{\cdots}{+} \\
	&\quad + [\iint_{\mathbf R^N \times \mathbf R^N}  \cdots \iint_{\mathbf R^N \times \mathbf R^N} ] [\varPsi(y,t) \lf( \slashed{\nabla}_{x,N-1} \slashed{\nabla}_{x,n} \delta(x-y) \rt)   S^\nabla_{N-1,n} (\psi) (x,t)]_{n=1}^N dxdy \\
	&\quad + [\iint_{\mathbf R^N \times \mathbf R^N}  \cdots \iint_{\mathbf R^N \times \mathbf R^N} ] [\varPsi(y,t) \lf( \slashed{\nabla}_{x,N} \slashed{\nabla}_{x,n} \delta(x-y) \rt)   S^\nabla_{N,n} (\psi) (x,t)]_{n=1}^N dxdy.
\end{align*}
On account of $\slashed{\nabla}_{x,1} \stackrel{\cdots}{+} \slashed{\nabla}_{x,N} = \sum_{n=1}^N \slashed{\nabla}_{x,n}$, we obtain
\begin{align*}
	(\textrm{Integral}_2) (t)
	&=\iint_{\mathbf R^N \times \mathbf R^N}  \varPsi(y,t) \sum_{n=1}^N \lf(\slashed{\nabla}_{x,n}^2 \delta(x-y) \rt) L^\nabla(\psi) (x,t)  dxdy  \\ 
	&\quad + \iint_{\mathbf R^N \times \mathbf R^N}  \varPsi(y,t) \sum_{n=1}^N \lf( \slashed{\nabla}_{x,1} \slashed{\nabla}_{x,n}\delta(x-y) \rt) S^\nabla_{1,n}(\psi) (x,t) dxdy \\ 
	&\quad + \iint_{\mathbf R^N \times \mathbf R^N} \varPsi(y,t) \sum_{n=1}^N \lf( \slashed{\nabla}_{x,2} \slashed{\nabla}_{x,n} \delta(x-y) \rt) S^\nabla_{2,n}(\psi) (x,t) dxdy \nn\\  
	&\quad \vdots \stackrel{\cdots}{+} \\ 
	&\quad+ \iint_{\mathbf R^N \times \mathbf R^N}  \varPsi(y,t)  \sum_{n=1}^N \lf( \slashed{\nabla}_{x,N-1} \slashed{\nabla}_{x,n} \delta(x-y)  \rt) S^\nabla_{N-1,n} (\psi) (x,t) dxdy  \\
	&\quad + \iint_{\mathbf R^N \times \mathbf R^N}  \varPsi(y,t) \sum_{n=1}^N \lf( \slashed{\nabla}_{x,N} \slashed{\nabla}_{x,n} \delta(x-y) \rt) S^\nabla_{N,n}(\psi) (x,t) dxdy .
\end{align*}
Using \eqref{lar} and summing up the above expressions we arrive at the following conclusion
\begin{align} \label{est:Int2}
	(\textrm{Integral}_2) (t)
	&= \iint_{\mathbf R^N \times \mathbf R^N}  \varPsi(y,t) \Delta_x \delta(x-y)  \lf(-\frac{1}{2} \mbox{div}\lf(\vec{\slashed{\nabla}}_x |\psi(x,t)|^2  \rt)  + \frac{\la \eta}{\eta+1} |\psi(x,t)|^{2\eta+2}  \rt) dxdy  \nn\\ 
	&\quad + \iint_{\mathbf R^N \times \mathbf R^N}  \varPsi(y,t) \sum_{m=1}^N \sum_{n=1}^N \lf(\slashed{\nabla}_{x,m} \slashed{\nabla}_{x,n}\delta(x-y) \rt) S^\nabla_{m,n}(\psi) (x,t) dxdy.
\end{align}
A combination \eqref{est:K}, \eqref{est:Int1} and \eqref{est:Int2} conduce us to the following estimate
\begin{align} \label{est:Int1+Int2}
	\frac{d}{dt} \mathcal V_\delta(t) &= (\textrm{Integral}_1)(t) + (\textrm{Integral}_2)(t) \nn \\
	&= \iint_{\mathbf R^N \times \mathbf R^N} \varPsi(y,t) \Delta_x \delta(x-y)  \lf(-\frac{1}{2} \mbox{div}\lf(\vec{\slashed\nabla}_x |\psi(x,t)|^2  \rt) \rt)   dxdy \nn\\
	&\quad + \iint_{\mathbf R^N \times \mathbf R^N} \varPsi(y,t) \Delta_x \delta(x-y)  \lf(\frac{\la \eta}{\eta+1} |\psi(x,t)|^{2\eta+2} \rt)   dxdy  \nn\\ 
	&\quad + \iint_{\mathbf R^N \times \mathbf R^N} \varPsi(y,t) \sum_{m=1}^N \sum_{n=1}^N \lf(\slashed{\nabla}_{x,m} \slashed{\nabla}_{x,n}\delta(x-y) \rt) S^\nabla_{m,n} (\psi)(x,t) dxdy \nn\\ 
	&\quad -2 \iint_{\mathbf R^N \times \mathbf R^N}  \lf(  \vec{\slashed{\nabla}}_{y} \cdot \overrightarrow{Q^\nabla} (\psi)(y,t) \rt) \lf(\vec{\slashed{\nabla}}_x \delta(x-y) \cdot \overrightarrow{Q^\nabla} (\psi)(x,t) \rt) dx  dy =  \sum_{k=3}^6 (\textrm{Integral}_k)(t) , 
\end{align}
where 
$$(\textrm{Integral}_3)(t) = \iint_{\mathbf R^N \times \mathbf R^N} \varPsi(y,t) \Delta_x \delta(x-y)  \lf(-\frac{1}{2} \mbox{div}\lf(\vec{\slashed\nabla}_x |\psi(x,t)|^2 \rt) \rt)   dxdy,$$ $$(\textrm{Integral}_4)(t) =  \iint_{\mathbf R^N \times \mathbf R^N} \varPsi(y,t) \Delta_x \delta(x-y)  \lf(\frac{\la \eta}{\eta+1} |\psi(x,t)|^{2\eta+2}\rt)   dxdy,$$
$$(\textrm{Integral}_5)(t) = \iint_{\mathbf R^N \times \mathbf R^N} \varPsi(y,t) \sum_{m=1}^N \sum_{n=1}^N \lf(\slashed{\nabla}_{x,m} \slashed{\nabla}_{x,n}\delta(x-y) \rt) S^\nabla_{m,n} (\psi)(x,t) dxdy,$$ $$(\textrm{Integral}_6)(t) = -2 \iint_{\mathbf R^N \times \mathbf R^N}  \lf(  \vec{\slashed{\nabla}}_{y} \cdot \overrightarrow{Q^\nabla} (\psi)(y,t) \rt) \lf(\vec{\slashed{\nabla}}_x \delta(x-y) \cdot \overrightarrow{Q^\nabla} (\psi)(x,t) \rt) dx  dy.$$
Note that the weight function $w$ here is a function of two variables $x,y$ and let us put $\delta(x-y) = |x-y|$, $x,y \in \mathbf R^N$, from the Appendix \ref{A1}, we have
$$\Delta_x |x-y| = \vec{\slashed\nabla}_x \cdot \vec{\slashed\nabla}_x |x-y| = \frac{N-1}{|x-y|}.$$
Thus, we have the estimate for the first integral in \eqref{est:Int1+Int2}, as follows
\begin{align*}
	(\textrm{Integral}_3)(t)
	&= \frac{N-1}{2} \iint_{\mathbf R^N \times \mathbf R^N} \frac{\varPsi(y,t)}{|x-y|} \lf(- \Delta_x \varPsi (x,t) \rt) dxdy \\
	&= \frac{N-1}{2} \int_{\mathbf R^N} (\int_{\mathbf R^N}\frac{\varPsi(y,t)}{|x-y|} dy) \lf(- \Delta_x \varPsi (x,t)\rt) dx.
\end{align*}
From the Appendix \ref{A2}, we derive that, for dimension $N\geq 1$
\begin{align} \label{est:Int3}
	(\textrm{Integral}_3)(t) &= \frac{N-1}{2} \int_{\mathbf R^N} (\int_{\mathbf R^N}\frac{\varPsi(y,t)}{|x-y|} dy) \lf(-\Delta_x \varPsi(x,t)\rt)  dx \nn\\
	&= \frac{N-1}{2} \int_{\mathbf R^N} ((-\Delta_x)^{-\frac{N-1}{2}} \varPsi(x,t)) (-\Delta_x \varPsi (x,t)) dx \nn\\
	&= \frac{N-1}{2} \int_{\mathbf R^N} ((-\Delta_x)^{-\frac{N-3}{2}} \varPsi^2(x,t)) dx \nn\\
	&= \frac{C_{N,\pi}(N-1)}{2} \lf\|(-\Delta_x)^{-\frac{N-3}{4}} \varPsi(\cdot,t) \rt\|_{L^2(\mathbf R^N)}^2 .
\end{align}
Now, we show that the sum of $(\textrm{Integral}_5)(t)$ and $(\textrm{Integral}_6)(t)$ is not negative. Indeed, by integrating by parts, we have
\begin{align*}
	(\textrm{Integral}_5)(t) + (\textrm{Integral}_6)(t) &= \iint_{\mathbf R^N \times \mathbf R^N} \varPsi(y,t) \sum_{m=1}^N \sum_{n=1}^N \lf( \slashed{\nabla}_{x,m} \slashed{\nabla}_{x,n}\delta(x-y) \rt) S^\nabla_{m,n} (\psi)(x,t) dxdy \\ 
	&\quad -2 \iint_{\mathbf R^N \times \mathbf R^N}  \sum_{m=1}^N \lf( \slashed{\nabla}_{y,m} Q^\nabla_m (\psi)(y,t) \rt) \sum_{n=1}^N \lf( \slashed{\nabla}_{x,n} \delta(x-y)  Q^\nabla_n (\psi) (x,t) \rt) dx  dy \\
	&= \iint_{\mathbf R^N \times \mathbf R^N}  \varPsi(y,t) \sum_{m=1}^N \sum_{n=1}^N \lf( \slashed{\nabla}_{x,m} \slashed{\nabla}_{x,n} \delta(x-y)   \rt) S^\nabla_{m,n} (\psi)(x,t) dxdy \\ 
	&\quad +2 \iint_{\mathbf R^N \times \mathbf R^N}  \sum_{m=1}^N \sum_{n=1}^N Q^\nabla_m (\psi)(y,t)  \lf( \slashed{\nabla}_{y,m} \slashed{\nabla}_{x,n} \delta(x-y)  Q^\nabla_n (\psi) (x,t)\rt) dx  dy .
\end{align*}
Thanks to Appendix \ref{A1}, one has  $\slashed{\nabla}_{y,m} \slashed{\nabla}_{x,n} |x-y| = - \slashed{\nabla}_{x,m} \slashed{\nabla}_{x,n} |x-y|$, for $x,y \in \mathbf R^N, m,n=1,...,N$, and we deduce that	
\begin{align*}
	&(\textrm{Integral}_5)(t) + (\textrm{Integral}_6)(t) \\
	&= 2\iint_{\mathbf R^N \times \mathbf R^N}   \sum_{m=1}^N \sum_{n=1}^N \lf( \slashed{\nabla}_{x,m} \slashed{\nabla}_{x,n}|x-y|  \rt) \\ &\quad\times \lf[\varPsi(y,t) \Re \lf(\slashed{\nabla}_{x,m} \psi (x,t) \slashed{\nabla}_{x,n} \bar{\psi} (x,t) \rt) - Q^\nabla_m (\psi)(y,t) Q^\nabla_n (\psi) (x,t)\rt] dx  dy \\
	&= 2\iint_{\mathbf R^N \times \mathbf R^N} \sum_{m=1}^N \sum_{n=1}^N \lf( \slashed{\nabla}_{x,m} \slashed{\nabla}_{x,n}|x-y| \rt) \\
	&\quad \times \lf[\frac{\varPsi(y,t)}{\varPsi (x,t)}\Re \lf(\bar{\psi}(x,t)\slashed{\nabla}_{x,m} \psi (x,t)  \psi(x,t) \slashed{\nabla}_{x,n} \bar{\psi} (x,t) \rt) - Q^\nabla_m (\psi)(y,t) Q^\nabla_n (\psi) (x,t)\rt] dx  dy.
\end{align*}	
By Appendix \ref{A1}, we have $\slashed{\nabla}_{x,m} \slashed{\nabla}_{x,n}|x-y| = \slashed{\nabla}_{y,m} \slashed{\nabla}_{y,n}|x-y|$, and because of the symmetry of the above expression, changing the roles of $x$ and $y$ does not change its value, so we deduce that
\begin{align} \label{est:I5+I6}
	&\quad (\textrm{Integral}_5)(t) + (\textrm{Integral}_6)(t) \nn\\
	&= \iint_{\mathbf R^N \times \mathbf R^N} \sum_{m=1}^N \sum_{n=1}^N \lf( \slashed{\nabla}_{x,m} \slashed{\nabla}_{x,n}|x-y|  \rt) \nn\\
	&\quad\times \lf[\frac{\varPsi(y,t)}{ \varPsi(x,t)} \Re \lf(\bar{\psi}(x,t)\slashed{\nabla}_{x,m} \psi (x,t)  \psi(x,t) \slashed{\nabla}_{x,n} \bar{\psi} (x,t) \rt) - Q^\nabla_m (\psi)(y,t) Q^\nabla_n (\psi) (x,t)\rt] dx  dy \nn\\
	&\quad+ \iint_{\mathbf R^N \times \mathbf R^N} \sum_{m=1}^N \sum_{n=1}^N \lf( \slashed{\nabla}_{y,m} \slashed{\nabla}_{y,n}|x-y|  \rt)  \nn\\
	&\quad\times \lf[\frac{\varPsi(x,t)}{\varPsi(y,t)}\Re \lf(\bar{\psi}(y,t) \slashed{\nabla}_{y,m} \psi (y,t) \psi(y,t)  \slashed{\nabla}_{y,n} \bar{\psi} (y,t) \rt) - Q^\nabla_m (\psi)(x,t) Q^\nabla_n (\psi) (y,t)\rt] dx  dy.
\end{align}
Putting $U= \bar{\psi}(x,t) \slashed{\nabla}_{x,m} \psi (x,t), V = \bar{\psi}(x,t) \slashed{\nabla}_{x,n} \psi (x,t)$, we have the following assertion
$$\Re(U\bar V) = \Re(U)\Re(V) + \Im(U)\Im(V),$$
and for complex number $z$, we also have $\Re(z\slashed{\nabla}_{x,m} \bar z) = \Re(\bar z\slashed{\nabla}_{x,m} z) = \frac{1}{2} \slashed{\nabla}_{x,m} |z|^2$, yields
\begin{align} \label{est:R1}
	\Re \lf(\bar{\psi}(x,t) \slashed{\nabla}_{x,m} \psi (x,t)  \psi(x,t)  \slashed{\nabla}_{x,n} \bar{\psi} (x,t) \rt) &=  \Re \lf(\bar{\psi}(x,t) \slashed{\nabla}_{x,m} \psi (x,t) \rt) \Re \lf(\bar{\psi}(x,t) \slashed{\nabla}_{x,n} \psi (x,t) \rt) \nn\\
	&\quad + \Im \lf(\bar{\psi}(x,t) \slashed{\nabla}_{x,m} \psi(x,t) \rt) \Im \lf(\bar{\psi}(x,t) \slashed{\nabla}_{x,n} \psi(x,t) \rt) \nn\\
	&= \frac{\slashed{\nabla}_{x,m} |\psi|^2(x,t)  \slashed{\nabla}_{x,n} |\psi|^2(x,t)}{4}  + Q^\nabla_m (\psi)(x,t) Q^\nabla_n (\psi)(x,t),
\end{align}
and the similarity to have
\begin{align} \label{est:R2}
	\Re \lf(\bar{\psi}(y,t) \slashed{\nabla}_{y,m} \psi (y,t)  \psi(y,t)  \slashed{\nabla}_{y,n} \bar{\psi} (y,t)  \rt) = \frac{\slashed{\nabla}_{y,m} |\psi|^2 (y,t) \slashed{\nabla}_{y,n} |\psi|^2 (y,t)}{4} + Q^\nabla_m (\psi)(y,t) Q^\nabla_n (\psi)(y,t).
\end{align}
From \eqref{est:I5+I6}, \eqref{est:R1} and \eqref{est:R2}, one obtains that
\begin{align*}
	&(\textrm{Integral}_5)(t) + (\textrm{Integral}_6)(t) = (\textrm{Integral}_7)(t) + (\textrm{Integral}_8)(t) + (\textrm{Integral}_9)(t) \\
	&=\iint_{\mathbf R^N \times \mathbf R^N} \sum_{m=1}^N \sum_{n=1}^N \lf( \slashed{\nabla}_{x,m} \slashed{\nabla}_{x,n} |x-y|  \rt) \frac{\varPsi(y,t)}{\varPsi(x,t)} \frac{\slashed{\nabla}_{x,m} |\psi|^2(x,t)  \slashed{\nabla}_{y,n} |\psi|^2(x,t) }{4} dx  dy \\
	&\quad  +  \iint_{\mathbf R^N \times \mathbf R^N} \sum_{m=1}^N \sum_{n=1}^N \lf( \slashed{\nabla}_{y,m} \slashed{\nabla}_{y,n} |x-y|  \rt) \frac{\varPsi(x,t)}{\varPsi(y,t)} \frac{\slashed{\nabla}_{y,m} |\psi|^2(y,t)  \slashed{\nabla}_{y,n} |\psi|^2(y,t) }{4}dx  dy \\
	&\quad + \iint_{\mathbf R^N \times \mathbf R^N} \sum_{m=1}^N \sum_{n=1}^N \lf( \slashed{\nabla}_{y,m} \slashed{\nabla}_{y,n}|x-y|  \rt)  \bigg[ \frac{\varPsi(y,t)}{\varPsi(x,t)} Q^\nabla_m (\psi)(x,t) Q^\nabla_n (\psi)(x,t) \\
	&\quad + \frac{\varPsi(x,t)}{\varPsi(y,t)} Q^\nabla_m (\psi)(y,t) Q^\nabla_n (\psi)(y,t)  -  Q^\nabla_m (\psi)(x,t) Q^\nabla_n (\psi)(y,t) - Q^\nabla_m (\psi)(y,t) Q^\nabla_n (\psi)(x,t) \bigg] dx  dy .
\end{align*}
Thanks to the Appendix \ref{A3} imply that the matrices $\slashed{\nabla}_{x,m} \slashed{\nabla}_{x,n}|x-y|$ and $\slashed{\nabla}_{y,m} \slashed{\nabla}_{y,n}|x-y|$ are positive definite, it follows that the integrals $(\textrm{Integral}_7)(t)$ and $(\textrm{Integral}_8)(t)$ are positive. Hence, 
\begin{align*}
	(\textrm{Integral}_5)(t) &+ (\textrm{Integral}_6)(t) 
	\geq 	\iint_{\mathbf R^N \times \mathbf R^N} \sum_{m=1}^N \sum_{n=1}^N \lf( \slashed{\nabla}_{y,m} \slashed{\nabla}_{y,n} |x-y|  \rt)  \bigg[\frac{\varPsi(y,t)}{\varPsi(x,t)} Q^\nabla_m (\psi)(x,t) Q^\nabla_n (\psi)(x,t)  \\
	& + \frac{\varPsi(x,t)}{\varPsi(y,t)} Q^\nabla_m (\psi)(y,t) Q^\nabla_n (\psi)(y,t) -  Q^\nabla_m (\psi)(x,t) Q^\nabla_n (\psi)(y,t) - Q^\nabla_m (\psi)(y,t) Q^\nabla_n (\psi)(x,t) \bigg] dx  dy .
\end{align*}
Defining the following vector function
$$\vec J(x,y,t) = \sqrt{\varPsi(y,t) \varPsi^{-1}(x,t)} \overrightarrow{Q^\nabla} (\psi)(x,t) - \sqrt{\varPsi(x,t) \varPsi^{-1}(y,t) } \overrightarrow{Q^\nabla} (\psi)(y,t),$$
we deduce that 
\begin{align} \label{est:Int5+Int6}
	(\textrm{Integral}_5)(t) + (\textrm{Integral}_6)(t)
	\geq \iint_{\mathbf R^N \times \mathbf R^N} \sum_{m=1}^N \sum_{n=1}^N \lf( \slashed{\nabla}_{y,m} \slashed{\nabla}_{y,n} |x-y|  \rt) J_m(x,y,t) J_n(x,y,t) dx  dy \geq 0.
\end{align}
From \eqref{est:Int1+Int2}, \eqref{est:Int3} and \eqref{est:Int5+Int6} it follows that
\begin{align*}
	\frac{d}{dt} \mathcal V_\delta(t) \geq \frac{C_{N,\pi}(N-1)}{2} \lf\|(-\Delta_x)^{-\frac{N-3}{4}} \varPsi(\cdot,t) \rt\|_{L^2(\mathbf R^N)}^2  + \frac{\la\eta(N-1)}{\eta+1} \iint_{\mathbf R^N \times \mathbf R^N} \frac{|\psi(y,t)|^2 |\psi(x,t)|^{2\eta+2} }{|x-y|}   dxdy  .
\end{align*}
Integrating from $0$ to $t$ and applying conservation of mass \eqref{ma}, we obtain
\begin{align*}	 
	\frac{C_{N,\pi}(N-1)}{2} \int_0^t \lf\|(-\Delta_x)^{-\frac{N-3}{4}} \varPsi(\cdot,z) \rt\|_{L^2(\mathbf R^N)}^2 dz  + \frac{\la\eta(N-1)}{\eta+1} \iiint_{[0,t] \times \mathbf R^N \times \mathbf R^N} \frac{|\psi(y,z)|^2 |\psi(x,z)|^{2\eta+2} }{|x-y|}   dxdydz \nn \\
	\leq \mathcal V_{|x-y|}(t) \leq \int_{\mathbf R^N} |\psi(y,t) |^2 dy \sup_{t \in \mathbf R^+ \atop y \in \mathbf R^N} \lf|\int_{\mathbf R^N} \vec{\slashed{\nabla}}_x |x-y| \cdot \overrightarrow{Q^\nabla} (\psi)(x,t) dx \rt| \nn \\
	= \|\psi(\cdot,t) \|_{L^2(\mathbf R^N)}^2 \sup_{t \in \mathbf R^+ \atop
		y \in \mathbf R^N} \lf| \int_{\mathbf R^N}  \vec{\slashed{\nabla}}_x |x-y| \cdot \overrightarrow{Q^\nabla} (\psi)(x,t) dx  \rt| \nn \\
	\stackrel{\eqref{ma}}{=} \|\psi_0\|_{L^2(\mathbf R^N)}^2 \sup_{t \in \mathbf R^+ \atop y \in \mathbf R^N} \lf|\int_{\mathbf R^N}  \vec{\slashed{\nabla}}_x |x-y| \cdot \overrightarrow{Q^\nabla} (\psi)(x,t) dx \rt|.
\end{align*}
The proof is complete. \hfill \rule{1.5mm}{3.5mm}	

\subsection{Proof of Proposition \ref{prop:space-time-est}} For $N\geq 1$, $\la>0$ and from \eqref{est:space-time}, we deduce that 
	\begin{align} \label{est:pro1}
	\lf\|(- \Delta_x)^{\frac{3-N}{4}} |\psi|^2 \rt\|_{L_t^2L_x^2(\mathbf R^+ \times \mathbf R^N)}  
	\lesssim  \lf\|\psi_0 \rt\|_{L^2(\mathbf R^N)} \sup_{t \in \mathbf R^+ \atop y \in \mathbf R^N} \sqrt{\lf|\int_{\mathbf R^N} \vec{\slashed{\nabla}}_x |x-y| \cdot \overrightarrow{Q^\nabla} (\psi)(x,t) dx \rt|}.
	\end{align}
From the Appendix \ref{A1}, we have that $\vec{\slashed{\nabla}}_x |x-y| = \frac{x-y}{|x-y|}$ and the fact that $\lf|\frac{x-y}{|x-y|}\rt| \leq 1$ for all $x,y \in \mathbf R^N$, yields
	$$\lf|\int_{\mathbf R^N} \vec{\slashed{\nabla}}_x |x-y| \cdot \overrightarrow{Q^\nabla} (\psi)(x,t) dx \rt| = \lf|\int_{\mathbf R^N} \frac{x-y}{|x-y|} \cdot \overrightarrow{Q^\nabla} (\psi)(x,t)  dx \rt| \leq \int_{\mathbf R^N} \lf| \overrightarrow{Q^\nabla} (\psi)(x,t) \rt| dx .$$	
By using the Cauchy-Schwartz inequality and from \eqref{def:Q}, we get that
	\begin{align} \label{est:pro2}
	\int_{\mathbf R^N} \lf| \overrightarrow{Q^\nabla} (\psi)(x,t) \rt| dx \leq \sqrt{ \sum_{m=1}^N  \lf| \int_{\mathbf R^N} \Im  (\bar\psi \slashed{\nabla}_{x,m} \psi)(x,t) dx\rt|^2  } \nn \\ \leq \sqrt{\int_{\mathbf R^N} \lf|\bar\psi(x,t)\rt|^2 dx  \sum_{m=1}^N  \int_{\mathbf R^N} \lf|\slashed{\nabla}_{x,m} \psi(x,t) \rt|^2 dx } \nn\\ 
	\leq \sqrt{\int_{\mathbf R^N} \lf|\psi(x,t)\rt|^2 dx}  \sqrt{\int_{\mathbf R^N} \lf|\vec{\slashed{\nabla}}_x \psi(x,t)\rt|^2 dx } = \lf\|\psi(\cdot,t) \rt\|_{L^2(\mathbf R^N)} \lf\|\vec{\slashed{\nabla}}_x \psi(\cdot,t) \rt\|_{L^2(\mathbf R^N)}.
	\end{align}
From \eqref{est:pro1}, \eqref{est:pro2} and conservation of mass \eqref{ma}, we obtain that
	\begin{align*}
	\lf\|(- \Delta_x)^{\frac{3-N}{4}} |\psi|^2\rt\|_{L_t^2L_x^2(\mathbf R^+ \times \mathbf R^N)}  
	\lesssim  \lf\|\psi_0 \rt\|_{L^2(\mathbf R^N)}^{\frac{3}{2}} \sup_{t \in \mathbf R^+} \lf\|\vec{\slashed{\nabla}}_x \psi \rt\|_{L^2(\mathbf R^N)} \\ \lesssim \lf\|\psi_0 \rt\|_{L^2(\mathbf R^N)}^{\frac{3}{2}} \sup_{t \in \mathbf R^+} \lf\|\vec{\slashed{\nabla}}_x \psi \rt\|_{L^2(\mathbf R^N)} \lesssim \lf\|\psi_0 \rt\|_{L^2(\mathbf R^N)}^{\frac{3}{2}} \sup_{t \in \mathbf R^+} \sqrt{\mathbb E[\psi](t)}.
	\end{align*}
From the conservation of energy \eqref{ener} one gets 
	\begin{align*}
	\lf\|(- \Delta_x)^{\frac{3-N}{4}} |\psi|^2\rt\|_{L_t^2L_x^2(\mathbf R^+ \times \mathbf R^N)}  \lesssim \lf\|\psi_0 \rt\|_{L^2(\mathbf R^N)}^{\frac{3}{2}}  \sqrt{\mathbb E[\psi_0]}.
	\end{align*}
The proof is done. \hfill \rule{1.5mm}{3.5mm}
\subsection{Proof of Theorem \ref{pseudo-conformal}} As argued above, Problem \eqref{Pro} possesses the solutions that satisfy the conservation laws \eqref{ener}, \eqref{mass} and \eqref{Mom}. The $(1-1)$ Morawetz action is the main tool to prove this theorem and we shall prove for $N=2$, the case $N=1$ is simpler and admitted. In fact, we recall that
	\begin{align} \label{def-Mo}
		\mathcal M_\delta(t) = \int_{\mathbf R^2} \vec \nabla \delta(x) \cdot \overrightarrow{Q^\nabla}(\psi)(x,t)  dx,
	\end{align}
	where the weight function $$\delta(\cdot): \mathbf R^2 \to \mathbf R \quad \mbox{is the moment arbitrary}.$$ 
By taking the derivative of both sides of \eqref{def-Mo} with respect to the time variable, we deduce
	\begin{align} \label{deri:M}
		\frac{d}{d t} \mathcal M_\delta(t) &= \int_{\mathbf R^2} \vec\nabla \delta(x) \cdot \frac{\partial \overrightarrow{Q^\nabla}(\psi)}{\partial t}(x,t)  dx \nn\\
		&= \int_{\mathbf R^2} \nabla_1 \delta(x) \frac{\partial Q^\nabla_1(\psi)}{\partial t}(x,t)    dx + \int_{\mathbf R^2} \nabla_2 \delta(x)  \frac{\partial Q^\nabla_2(\psi)}{\partial t} (x,t) dx  .
	\end{align}
From the local momentum conservation \eqref{Mom}, we have for $m = 1,2$
	\begin{align} \label{deri:J}
		\frac{\partial Q^\nabla_m(\psi)}{\partial t} (x,t) &= -  \nabla_1 \lf(\mathbf 1_{m, 1} L^\nabla(\psi)  + S^\nabla_{m, 1} (\psi) \rt)(x,t)  - \nabla_2 \lf(\mathbf 1_{m, 2} L^\nabla(\psi)  + S^\nabla_{m, 2} (\psi) \rt)(x,t) .
	\end{align}
From the definition of Kronecker delta function \eqref{Kronecker} and combining \eqref{deri:M} and \eqref{deri:J}, we obtain that
	\begin{align*}
	\frac{d}{d t} \mathcal M_\delta(t)
	&= - \int_{\mathbf R^2} \nabla_1 \delta(x)  \nabla_1 \lf(L^\nabla(\psi) +  S^\nabla_{1,1}(\psi) \rt)(x,t)   dx - \int_{\mathbf R^2} \nabla_1 \delta(x)  \nabla_2 S^\nabla_{1,2} (\psi)(x,t)  dx  \nn\\
	&\quad - \int_{\mathbf R^2} \nabla_2 \delta(x)  \nabla_1 S^\nabla_{2,1} (\psi)(x,t)   dx - \int_{\mathbf R^2} \nabla_2 \delta(x)  \nabla_2 \lf(L^\nabla(\psi) +  S^\nabla_{2,2}(\psi) \rt)(x,t) dx.
	\end{align*}
From this, integrating by parts, we derive
	\begin{align*}
	\frac{d}{d t} \mathcal M_\delta(t) &= \int_{\mathbf R^2} \nabla_1^2  \delta(x) L^\nabla(\psi)(x,t) dx + \int_{\mathbf R^2}  \nabla_2^2  \delta(x) L^\nabla(\psi)(x,t) dx \nn\\
	&\quad + \int_{\mathbf R^2} \nabla_1^2 \delta(x) S^\nabla_{1,1} (\psi)(x,t) dx + \int_{\mathbf R^2}  \nabla_1 \nabla_2 \delta(x) S^\nabla_{1,2} (\psi)(x,t) dx \nn\\
	&\quad + \int_{\mathbf R^2} \nabla_2 \nabla_1 \delta(x) S^\nabla_{2,1} (\psi)(x,t) dx + \int_{\mathbf R^2} \nabla_2^2\delta(x) S^\nabla_{2,2} (\psi)(x,t) dx .
	\end{align*}
By applying the definition of the Kronecker delta function \eqref{Kronecker}, we conclude that
	\begin{align*}
		\frac{d}{d t} \mathcal M_\delta(t) &=  \int_{\mathbf R^2} \lf(\nabla_1^2 \delta(x) + \nabla_2^2 \delta(x) \rt) L^\nabla(\psi)(x,t) dx \nn\\ 
		&\quad + \int_{\mathbf R^2} \sum_{n=1,2}  \nabla_1 \nabla_n \delta(x)   S^\nabla_{1,n}(\psi) (x,t) dx + \int_{\mathbf R^2} \sum_{n=1,2}  \nabla_2 \nabla_n \delta(x)  S^\nabla_{2,n}(\psi)(x,t) dx  .
	\end{align*}
	Consequently,
	\begin{align} \label{dM}
		\frac{d}{d t} \mathcal M_\delta(t) =  \int_{\mathbf R^2} \mbox{div} \lf(\vec\nabla \delta(x)\rt) L^\nabla(\psi)(x,t)  dx + \int_{\mathbf R^2}  \sum_{m=1,2} \sum_{n=1,2} \nabla_m \nabla_n \delta(x) S^\nabla_{m,n} (\psi)(x,t)  dx.
	\end{align}
	Introducing \eqref{lar} and \eqref{ten} into \eqref{dM},  we obtain that
	\begin{align} \label{dMa}
		\frac{d}{d t} \mathcal M_\delta(t) &=  \int_{\mathbf R^2} \mbox{div} \lf(\vec \nabla \delta(x) \rt) \lf(-\frac{1}{2} \mbox{div} \lf(\vec \nabla |\psi(x,t)|^2 \rt) +  \frac{\la \eta}{\eta+1} |\psi(x,t)|^{2\eta+2}  \rt)  dx \nn\\ 
		&\quad + \int_{\mathbf R^2} \sum_{m =1,2} \sum_{n=1,2} \lf( \nabla_m \nabla_n \delta(x) \rt) \lf(2 \Re  \lf(\nabla_m \psi \nabla_n \bar{\psi} \rt)(x,t) \rt)   dx.
	\end{align}
	Now, we pick the weight function $\delta(x) = |x|^2$, from the Appendix \ref{A0} we have  that
	$$\mbox{div} \lf(\vec\nabla \delta \rt) = \Delta |x|^2 = 4, \quad \mbox{and}\quad \nabla_m \nabla_n |x|^2 = 2 \times \mathbf 1_{m, n}, \quad m = 1,2, n = 1,2. $$
	Combining this comment with  inequality \eqref{dMa} and the fact that $$\sum_{m =1,2} \sum_{n=1,2} \lf( \nabla_m \nabla_n |x|^2\rt) \lf(2 \Re  \lf(\nabla_m \psi \nabla_n \bar{\psi} \rt)(x,t) \rt) = \sum_{m =1,2} \sum_{n=1,2}  2 (\mathbf 1_{m, n}) \lf(2 \Re  \lf(\nabla_m \psi \nabla_n \bar{\psi} \rt) (x,t) \rt) = 4 |\vec \nabla \psi(x,t)|^2,$$ we deduce
	\begin{align*} 
		\frac{d}{d t} \mathcal M_{\delta = |x|^2}(t) =  \int_{\mathbf R^2} \frac{4\la \eta}{\eta+1} |\psi(x,t)|^{2\eta+2}  dx - 2 \int_{\mathbf R^2}  \mbox{div} \lf(\vec \nabla |\psi(x,t)|^2  \rt)  dx 
		+ 4\int_{\mathbf R^2} |\vec \nabla \psi (x,t) |^2 dx .
	\end{align*}
	Equivalent to
	\begin{align*}
		\frac{d}{d t} \mathcal M_{\delta = |x|^2}(t) 
		&= 8 \lf(\frac{1}{2} \int_{\mathbf R^2}  |\vec \nabla \psi  (x,t)|^2 dx + \frac{\la}{2\eta+2} \int_{\mathbf R^2} |\psi(x,t)|^{2\eta+2} dx\rt)  \nn\\ 
		&\quad - \frac{\la (4 - 4 \eta)}{\eta+1}  \int_{\mathbf R^2} |\psi(x,t)|^{2\eta+2}  dx - 2 \int_{\mathbf R^2} \mbox{div} \lf(\vec \nabla |\psi(x,t)|^2 \rt) dx .
	\end{align*}
	By the total energy $\mathbb E[\psi] (t)$ in \eqref{ener}, one obtains 
	\begin{align} \label{Mo-deri}
		\frac{d}{d t} \mathcal M_{\delta = |x|^2}(t) = 8 \mathbb E[\psi] (t) - \frac{\la(4 - 4\eta) }{\eta+1}  \int_{\mathbf R^2}  |\psi(x,t)|^{2\eta+2}  dx  
		- 2 \int_{\mathbf R^2}  \mbox{div} \lf(\vec \nabla |\psi(x,t)|^2 \rt) dx . 
	\end{align}
Let us define the $C^2$ function by the quantity   
	$$\mathfrak S(t) =  \int_{\mathbf R^2} |x|^2 \varPsi(x,t) dx, \quad\mbox{with}~ \varPsi = |\psi|^2.$$
From \eqref{mass}, we have that
	$$\frac{d}{d t} \mathfrak S (t) = \int_{\mathbf R^2}  |x|^2 \frac{\partial \varPsi }{\partial t}(x,t) dx = -2 \int_{\mathbf R^2}  |x|^2 \lf(\vec \nabla \cdot \overrightarrow{Q^\nabla}(\psi) (x,t)\rt)dx = 2 \mathcal M_{\delta = |x|^2}(t),$$	
where we have used the integration by parts. Thus we conclude that 
	\begin{align*}
		\frac{d^2}{d t^2} \mathfrak S(t) = 16 \mathbb E[\psi] (t) - \frac{8\la(1 - \eta)}{\eta +1}  \int_{\mathbf R^2} |\psi(x,t)|^{2\eta+2}  dx  
		- 4 \int_{\mathbf R^2}  \mbox{div} \lf(\vec \nabla |\psi(x,t)|^2 \rt) dx . 
	\end{align*}
From \eqref{def-Pseudo}, we recall that
	\begin{align*}
		\mathbb{P}[\psi](t) = \displaystyle \lf\|(\cdot + 2\mathrm it \vec\nabla) \psi (\cdot,t)\rt\|^2_{L^2(\mathbf R^2)}   + \frac{4t^2\la}{\eta+1} \int_{\mathbf R^2}  |\psi(x,t)|^{2\eta+2} dx  
		- 4 \int_{\mathbf R^2} \mbox{div} \lf(\vec \nabla \mathcal H(\psi) (x,t)\rt) dx .
	\end{align*}
The above expression is equivalent to
	\begin{align*} 
		\mathbb{P}[\psi](t) &= \lf\|(\cdot) \psi (\cdot,t) \rt\|^2_{L^2(\mathbf R^2)} + 4t^2 \lf\|\vec \nabla \psi (\cdot,t) \rt\|^2_{L^2(\mathbf R^2)} - 4t \int_{\mathbf R^2} x \cdot \overrightarrow{Q^\nabla}(\psi)(x,t) dx  \nn \\ 
		&\quad + \frac{4t^2\la}{\eta+1} \int_{\mathbf R^2} |\psi(x,t)|^{2\eta+2} dx - 4 \int_{\mathbf R^2}  \mbox{div} \lf(\vec \nabla \mathcal H(\psi) (x,t)\rt) dx   .  
	\end{align*}
If we understand that $|x|^2 = \delta(x)$ and using the  conservation of energy \eqref{ener} then the equation above is
	\begin{align} \label{def-P}
		\mathbb{P}[\psi](t) & \stackrel{\eqref{ener}}{=} \int_{\mathbf R^2}  \delta(x)|\psi(x,t)|^2 dx + 8t^2 \mathbb E[\psi_0]  - 2 t \int_{\mathbf R^2} \vec \nabla \delta(x) \cdot \overrightarrow{Q^\nabla}(\psi)(x,t) dx  - 4 \int_{\mathbf R^2} \mbox{div} \lf(\vec \nabla \mathcal H(\psi) (x,t) \rt) dx .  
	\end{align}
The derivative of the above equation according to the variable $t$ becomes
	\begin{align*}
		\frac{d}{d t} \mathbb{P}[\psi] (t) 
		&= \frac{d}{dt} \int_{\mathbf R^2} \delta(x) \varPsi(x,t) dx + 16 t \mathbb E[\psi_0] \nn\\
		&\quad - 2 \frac{d}{dt}\lf( t \int_{\mathbf R^2} \vec \nabla \delta(x) \cdot \overrightarrow{Q^\nabla}(\psi)(x,t) dx \rt)  - 4  \int_{\mathbf R^2} \mbox{div} \lf(\vec \nabla \lf( \frac{\partial}{\partial t} \mathcal H(\psi) (x,t) \rt)\rt) dx, \quad \mbox{for}~ \varPsi  = |\psi|^2.
	\end{align*}
Moreover, from \eqref{mass} and using integration by parts, we have
	$$\frac{d}{d t} \int_{\mathbf R^2}  \delta \varPsi dx =  \int_{\mathbf R^2} \delta \frac{\partial  \varPsi}{\partial t}  dx =   2 \int_{\mathbf R^2} \vec \nabla \delta \cdot \overrightarrow{Q^\nabla}(\psi)  dx,$$
and \begin{align*}
		\frac{d}{dt} \lf(t \int_{\mathbf R^2} \vec\nabla \delta \cdot \overrightarrow{Q^\nabla}(\psi)  dx \rt)  = t \int_{\mathbf R^2} \vec \nabla \delta \cdot \frac{\partial}{\partial t} \overrightarrow{Q^\nabla}(\psi)  dx +   \int_{\mathbf R^2} \vec \nabla \delta \cdot \overrightarrow{Q^\nabla}(\psi)  dx.
	\end{align*}
Consequently, 
	\begin{align*}
		\frac{d}{d t} \mathbb{P}[\psi] (t) 
		&= -2t \int_{\mathbf R^2} \vec\nabla \delta(x) \cdot \frac{\partial}{\partial t} \overrightarrow{Q^\nabla}(\psi)(x,t)  dx  + 16t \mathbb E[\psi_0]  - 4 \int_{\mathbf R^2} \mbox{div} \lf(\vec\nabla \lf(\frac{\partial}{\partial t} \mathcal H(\psi) (x,t)\rt) \rt) dx.
	\end{align*}
From the above observation and combined with \eqref{def-Mo}, \eqref{Mo-deri}, recalling that $\delta(x) = |x|^2$, one obtains 
	\begin{align} \label{P}
		\frac{d}{d t} \mathbb{P}[\psi] (t) 
		&= -2t \frac{d}{d t} \mathcal M_{\delta = |x|^2}(t) + 16t \mathbb E[\psi_0] - 4 t \int_{\mathbf R^2} \mbox{div} \lf(\vec \nabla |\psi(x,t)|^2 \rt) dx \nn\\
		& = \frac{2\la t}{\eta+1}(4 - 4\eta) \int_{\mathbf R^2}  |\psi(x,t)| ^{2\eta+2} dx ,
	\end{align}
	where we have used the conditions \eqref{cond:H} that $\frac{\partial \mathcal H(\psi)}{\partial t}  (x,t) = t |\psi|^2(x,t)$. Let us pick $\eta= 1$ in \eqref{P} and integrating this equation in time on $[0,t]$ and we obtain that  $\mathbb{P}[\psi](t)$ is conserved. The proof of  \autoref{pseudo-conformal} is finished.
	\hfill \rule{1.5mm}{3.5mm}
	\begin{remark}
		Obviously, we can choose the weight function $\delta(\cdot): \mathbf R^N \to \mathbf R$ in other forms. However, in this proof, we choose $$\delta =|x|^2 = \sum_{k=1}^N x_k^2,$$ because of the properties for our convenience such as:
		$$\mathrm{div}\lf(\nabla |x|^2\rt) =\frac{\partial^2 |x|^2}{\partial x_1^2} + \frac{\partial^2 |x|^2}{\partial x_2^2} \stackrel{\cdots}{+} \frac{\partial^2 |x|^2}{\partial x_N^2}  = 2N,$$ and $$\frac{\partial }{\partial x_m} \frac{\partial }{\partial x_n} |x|^2 =2 \times \mathbf 1_{m n}, \quad \mbox{for}~~ 1 \leq m ,n \leq N.$$
	\end{remark}
	\subsection{Proof of Theorem \ref{Morawetz estimate}}
	The proof of this theorem can be seen as a consequence of \autoref{pseudo-conformal}. Indeed, from \eqref{P},  we have that
	\begin{align*} 
		\frac{2\la t}{\eta+1}(4 - 2\eta ) \int_{-\infty}^\infty |\psi(x,t)|^{2\eta+2}  dx = \frac{d}{d t} \mathbb{P}[\psi] (t).
	\end{align*}
	For $\eta  < 2$, for all $t \in [0,T]$, we have
	\begin{align} \label{est:space}
		\int_{-\infty}^\infty  t |\psi(x,t)|^{2\eta+2} dx   \lesssim \frac{d}{d t} \mathbb{P}[\psi] (t).
	\end{align}
Integrating this equation in time on $[0,T]$ we obtain
	\begin{align*}
	\int_0^T \int_{-\infty}^{\infty}  t |\psi(x,t)|^{2\eta+2}  dx dt  &\lesssim \mathbb{P}[\psi] (t) - \mathbb{P}[\psi] (0) \\
 	& \lesssim \mathbb{P}[\psi] (t) - \lf\|(\cdot) \psi_0\rt\|^2_{L^2(\mathbf R)}. 
 	\end{align*} 
Applying the fundamental theorem of calculus we deduce \eqref{Est}. The proof of \autoref{Morawetz estimate} is complete.
	\hfill \rule{1.5mm}{3.5mm}
	
	\subsection{Proof of Proposition \ref{decay}} This proof is also derived from the proof of \autoref{pseudo-conformal}. In fact,
	for $\eta \in [1,2)$, for all $t \in (0,T]$, from \eqref{est:space} we have
	$$ \int_{-\infty}^\infty  t |\psi(x,t)|^{2\eta +2} dx   \lesssim  \sup_{0<t \leq T} \lf(\frac{d}{d t} \mathbb{P}[\psi] (t) \rt).$$
	Assume that $\displaystyle\sup_{0<t\leq T} \lf(\frac{d}{d t} \mathbb{P}[\psi] (t) \rt)$ has an upper bound $C$ for some $C>0$, then 
	$$ \int_{-\infty}^\infty t |\psi(x,t)|^{2\eta +2} dx  \leq C .$$
	Consequently, 
	$$ t^{\frac{1}{2\eta +2}} \lf\|\psi(\cdot,t)\rt\|_{L^{2\eta +2}(\mathbf R)} \leq C .$$
	This finishes the proof of \autoref{decay}.
	\hfill \rule{1.5mm}{3.5mm}
	
	\section{Conclusions}
	This study has achieved a major extension over the much more investigated nonlinear Schr\"odinger equation. Using the Gagliardo-Nirenberg interpolation inequality, we obtained the local well-posedness results of the solutions: local existence, continuation, blow-up alternative and continuous dependence on initial data. Moreover, by using the two types of Morawetz actions, we established the stronger upper bound of the solutions, and moreover the modified pseudo-conformal conservation law for Problem \eqref{Pro}. The Morawetz and decay estimates are also investigated.\\
	
	\noindent \textbf{\large{Acknowledgment}.}
	V.V. Au and N.H. Tuan are supported by the Van Lang University.

\section{Appendix}
\subsection{A1} \label{A0} On the $N$-dimensional Euclidean space $\mathbf R^N$, let $x = (x_1,...,x_N) \in \mathbf R^N$, we define the Euclidean norm of $x$ on the $N$-dimensional vector spaces with
	$|x| = \sqrt{\sum_{k=1}^N x_k^2}.$
Then we have $\vec\nabla |x|^2 = \lf(\nabla_1|x|^2,...,\nabla_N|x|^2\rt)=(2x_1,...,2x_N)$ and 
	$\Delta |x|^2 = (\vec\nabla \cdot \vec\nabla) |x|^2 = 2N,$ and
	$$\nabla_m \nabla_n |x|^2 = \begin{cases}
		2, \quad \mbox{if}\quad m=n \\
		0, \quad \mbox{if} \quad m \neq n
	\end{cases} =2 \times \mathbf 1_{m,n}, \quad 1 \leq m,n \leq N.$$
\subsection{A2} \label{A1}	 For $x=(x_1,...,x_N), y = (y_1,...,y_N) \in \mathbf R^N$, we have the Euclidean norm of $x-y$ as
$$|x-y| = \sqrt{\sum_{k=1}^N (x_k-y_k)^2},$$
and $\vec{\slashed{\nabla}}_x |x-y| = (\slashed{\nabla}_{x,1} |x-y|, ..., \slashed{\nabla}_{x,N} |x-y|)$, then we deduce 
$$\vec{\slashed{\nabla}}_x |x-y| = \lf(\frac{x_1-y_1}{\sqrt{\sum_{k=1}^N (x_k-y_k)^2}},...,\frac{x_N-y_N}{\sqrt{\sum_{k=1}^N (x_k-y_k)^2}}\rt) ,$$
and we have
\begin{align*}
\Delta_x |x-y| = (\vec{\slashed{\nabla}}_x \cdot \vec{\slashed{\nabla}}_x) |x-y| = \sum_{k=1}^N \frac{\sum_{j=1}^N (x_j-y_j)^2 - (x_k-y_k)^2}{\lf(\sqrt{\sum_{j=1}^N (x_j-y_j)^2}\rt)^3}= \frac{N-1}{\sqrt{\sum_{k=1}^N (x_k-y_k)^2}}.
\end{align*}
We can also check that 
$$\slashed{\nabla}_{y,m} \slashed{\nabla}_{x,n} |x-y| = - \slashed{\nabla}_{x,m} \slashed{\nabla}_{x,n} |x-y|, $$ and $$\slashed{\nabla}_{x,m} \slashed{\nabla}_{x,n} |x-y| = \slashed{\nabla}_{y,m} \slashed{\nabla}_{y,n} |x-y|.$$
\subsection{A3} \label{A3}  We shall prove that the following matrix is positive definite
	$$\lf(\begin{array}{cccc}
	\slashed\nabla_{x,1}^2 |x-y|& \slashed\nabla_{x,1}\slashed\nabla_{x,2} |x-y| &\dots&\slashed\nabla_{x,1}\slashed\nabla_{x,N} |x-y|  \\
	\slashed\nabla_{x,2}\slashed\nabla_{x,1} |x-y|& \slashed\nabla_{x,2}^2 |x-y|&\dots&\slashed\nabla_{x,2}\slashed\nabla_{x,N} |x-y|\\
	\vdots & \vdots&\ddots&\vdots\\
	\slashed\nabla_{x,N}\slashed\nabla_{x,1} |x-y|&\slashed\nabla_{x,N}\slashed\nabla_{x,2} |x-y|&\dots&\slashed\nabla_{x,N}^2 |x-y|
	\end{array}	\rt)_{N\times N}$$
Picking the vector function $\{\vartheta_k\}_{k=1}^N$, with $\vartheta_k (\cdot): \mathbf R^N \to \mathbf C$,  we need to prove that
	\begin{align} \label{A3-1}
	\int_{\mathbf R^N} \sum_{m=1}^N\sum_{n=1}^N \lf(\slashed\nabla_{x,m}\slashed\nabla_{x,n} |x-y|\rt) \vartheta_m(x)\vartheta_n(x) dx \geq 0.
	\end{align}
We have that $$\slashed\nabla_{x,m}\slashed\nabla_{x,n} |x-y| = \frac{1}{|x-y|} \lf(\mathbf 1_{mn} - \frac{(x_m-y_m)(x_n-y_n)}{|x-y|^2}\rt).$$
Now, we can write $x,y$ as the vector functions in $\mathbf R^N$ by $\vec x = (x_1,...,x_N), \vec y = (y_1,...,y_N)$ and $\vec \vartheta = (\vartheta_1,...,\vartheta_N)$ and one obtains
	\begin{align} \label{A3-2}
	&\sum_{m=1}^N\sum_{n=1}^N \lf(\slashed\nabla_{x,m}\slashed\nabla_{x,n} |\vec x- \vec y|\rt) \vartheta_m(x)\vartheta_n(x) \nn\\
	&= \sum_{m=1}^N\sum_{n=1}^N \frac{1}{|\vec x-\vec y|} \lf(\mathbf 1_{mn} - \frac{(x_m-y_m)(x_n-y_n)}{|\vec x-\vec y|^2}\rt) \vartheta_m(x)\vartheta_n(x) \nn\\
	&=\sum_{m=1}^N\sum_{n=1}^N \frac{1}{|\vec x-\vec y|} \lf(\lf|\vec \vartheta (x)\rt|^2 - \frac{\lf[(\vec{x } - \vec{y}) \cdot \vec{\vartheta} \rt]^2}{|\vec x-\vec y|^2}\rt).
	\end{align}
From the Cauchy-Schwartz inequality, we have 
	\begin{align} \label{A3-3}
	\frac{\lf[(\vec{x }- \vec{y}) \cdot \vec{\vartheta} \rt]^2}{|\vec x-\vec y|^2} &= \frac{\lf[(x_1-y_1) \vartheta_1 \stackrel{\cdots}{+} (x_N-y_N)\vartheta_N \rt]^2}{|\vec x-\vec y|^2} \nn\\
	&\leq \frac{\lf[(x_1-y_1)^2 \stackrel{\cdots}{+} (x_N-y_N)^2\rt](\vartheta_1^2 \stackrel{\cdots}{+} \vartheta_N^2)}{|\vec x-\vec y|^2} \leq \lf|\vec \vartheta\rt|^2.
	\end{align}
Thus \eqref{A3-1} follows from \eqref{A3-2} and \eqref{A3-3}. Similarly, we infer that the matrix $\lf\{\slashed\nabla_{y,m}\slashed\nabla_{y,n} |x-y| \rt\}_{m,n=1}^N$ is also positive definite.    
\subsection{A4} \label{A2} Let $f \in L^1(\mathbf R^N)$, we have the Fourier transform of $f$ as
	$$\widehat f (\xi) = \int_{\mathbf R^N} f(x) e^{-2\pi \mathrm i \sum_{k=1}^N x_k\xi_k } dx.$$
By integrating by parts, we have 
	$$(\widehat{-\Delta_x f}) (\xi) = - \int_{\mathbf R^N} \Delta_x f (x) e^{-2\pi \mathrm i \xi \cdot x} dx = C_\pi |\xi|^2 \widehat f (\xi).$$
We also have that $\widehat{|x|^{-1}} = C_{N,\pi} |\xi|^{-(N-1)}$, $N>1$, and with $C_{N,\pi} = 4\pi^{\frac{3}{2}} \frac{\Gamma(\frac{1}{2})}{\Gamma(\frac{N-1}{2})}$. Then we may define the identity
	$$(-\Delta_x)^{-\frac{N-1}{2}} f (x)  = C_{N,\pi} \int_{\mathbf R^N} \frac{f(y)}{|x-y|} dy.$$
And if $f \in C^\infty_c(\mathbf R^N)$ then  $f \in L^1(\mathbf R^N) \cap L^2(\mathbf R^N)$ and
	$$\|\widehat f\|_{L^2(\mathbf R^N)} = C_{N,\pi} \| f\|_{L^2(\mathbf R^N)}, \quad C_{N,\pi} = (2\pi)^N.$$


\begin{thebibliography}{99}
		\bibitem{AM013}  Ablowitz, M.J.,  Musslimani, Z.H.: Integrable nonlocal nonlinear Schr\"odinger equation. Phys. Rev. Lett. \textbf{110}(6) (2013), 064105.	
			
		\bibitem{B004}  Banica, V.: Remarks on the blow-up for the Schr\"odinger equation with critical mass on a plane domain. Ann. Scuola Norm. Sup. Pisa Cl. Sci. \textbf{5}(1) (2004), 139--170.
		
		\bibitem{BR015}  Bisci,  G.M.,  R\u{a}dulescu, V.D.: Ground state solutions of scalar field fractional Schr\"odinger equations. Calc. Var. Partial Differ. Equ. \textbf{54} (2015), 2985--3008.
		
		\bibitem{Bo99}  Bourgain, J.: Global wellposedness of defocusing critical nonlinear Schr\"odinger equation in the radial case.  J.
			Amer. Math. Soc.  \textbf{12}(1) (1999), 145--171.
		
		\bibitem{BM018}  Br\'{e}zis, H.,  Mironescu, P.: Gagliardo-Nirenberg inequalities and non-inequalities: The full story.  Ann. I. H. Poincar\'{e} - AN.  \textbf{35} (2018), 1355--1376.
		
		\bibitem{BGT003}  Burq, N.,  G\'{e}rard, P.,  Tzvetkov, N.: Two singular dynamics of the nonlinear Schr\"odinger
		equation on a plane domain.  Geom. Funct. Anal. \textbf{13}(1) (2003), 1--19.
		
		\bibitem{C003}  Cazenave, T.: Semilinear Schr\"odinger equation. Courant lecture notes 10.  Amer. Math. Soc. (2003), 346 pages.
		
		\bibitem{CN021}  Campos, B.J.,  Naumkin, P.I.: Large time asymptotics for the higher-order nonlinear nonlocal Schr\"odinger equation.  Nonlinear Anal.
		\textbf{205} (2021), 112238.
		
		\bibitem{CTG009}  Colliander, J.,   Tzirakis, N.,  Grillakis, M.G.: Tensor products and correlation estimates with applications to nonlinear Schr\"odinger equations.  Commun. Pure Appl. Math.  \textbf{62}(7) (2009), 920--968.
		
		\bibitem{CRSW004} Colliander, J.,  Raynor, S., Sulem, C.,  Wright, J.D.: Ground state mass concentration in the $L^2$-critical nonlinear Schr\"odinger equation below $H^1$.  Math. Research Lett.  \textbf{12}(3) (2004).
		
		\bibitem{CGT010}  Colliander, J.,  Grillakis, M.G., Tzirakis, N.: Remarks on global a priori estimates for the Schr\"odinger equation.  Proc. Amer. Math. Soc. 
		\textbf{138}(12) (2010), 4359--4371.
		
		\bibitem{CKSTT008}  Colliander, J.,  Keel, M.,  Staffilani, G.,  Takaoka, H.,  Tao, T.: Global well-posedness and scattering
		for the energy-critical nonlinear
		Schr\"odinger equation in $\mathbf R^3$. Ann. Math.  \textbf{167} (2008), 767--865. 
		
		\bibitem{DW018} Ding, Q., Wang, Y.: Vortex filament on symmetric Lie algebras and generalized bi-Schr\"odinger flows. Math. Z., \textbf{290} (2018), 167--193.
		
		\bibitem{D019} Duong, D.V.: Energy scattering for a class of the defocusing inhomogeneous nonlinear Schr\"odinger equation. J. Evol. Equ. \textbf{19} (2019), 411--434.
		
		\bibitem{D018} Duong, D.V.: Blowup of $H^1$ solutions for a class of the focusing inhomogeneous nonlinear Schr\"odinger equation. Nonlinear Anal. \textbf{174} (2018), 169--188.
		
		\bibitem{FW021} Fern\'{a}ndez, A.J., Weth, T.: The nonlinear Schr\"odinger equation in the half-space. Math. Ann. (2021)
		
		\bibitem{GV79-I} Ginibre, J., Velo, G.: On a class of nonlinear Schr\"odinger equations. I. The Cauchy problem, general
		case. J. Funct. Anal. \textbf{32}(1) (1979), 1--32.
		
		\bibitem{GV79-II}  Ginibre, J.,  Velo, G.: On a class of nonlinear Schr\"odinger equations.
		II. Scattering theory, general case. J. Funct. Anal. \textbf{32}(1) (1979), 33--71.
		
		\bibitem{GV80} Ginibre, J., Velo, G.: On a class of non linear Schr\"odinger equations with non local interaction. Math. Z. \textbf{170}, 109–136 (1980)
		
		\bibitem{JL021} Jeanjean, L., Le, T.T.: Multiple normalized solutions for a Sobolev critical Schr\"odinger equation. Math. Ann. (2021).
		
		\bibitem{KOPV017}  Killip, R.,  Oh, T.,  Pocovnicu, O., Visan, M.: Solitons and scattering for the cubic–quintic nonlinear Schr\"odinger equation on $\mathbf R^3$. Arch. Ration. Mech. Anal.  \textbf{225} (2017), 469--548.
		
		
		\bibitem{KV010}  Killip, R.,  Visan, M.: The focusing energy-critical nonlinear Schr\"odinger equation in dimensions five and higher. Amer. J. Math. \textbf{132}(2) (2010), 361--424.
		
		\bibitem{KMVZZ018}  Killip, R.,  Miao, C.,  Visan, M.,  Zhang, J.,  Zheng, J.: Sobolev spaces adapted to the Schr\"odinger operator with inverse-square potential. Math. Z. \textbf{288} (2018), 1273--1298.
		
		\bibitem{KMV018}  Killip, R.,  Murphy, J.,   Visan, M.: The initial-value problem for the cubic-quintic NLS with nonvanishing boundary conditions. SIAM J. Math. Anal.  \textbf{50}(3) (2018), 2681--2739.
		
		\bibitem{MS021} Martin, J., Pravda-Starov, K.: Geometric conditions for the exact controllability of fractional free and harmonic Schr\"odinger equations. J. Evol. Equ. \textbf{21} (2021), 1059--1087.
		
		\bibitem{MR008}  Merle, F.,  Rapha\"{e}l, P.: Blow up of the critical norm for some radial $L^2$ super critical nonlinear Schr\"odinger equations. Amer J. Math. 
		\textbf{130}(4) (2008), 945--978.
		
		\bibitem{OW020}  Oh, T.,  Wang, Y.: Global well-posedness of the one-dimensional cubic nonlinear Schr\"odinger equation in almost critical spaces. J. Differ. Equ.
		\textbf{269}(1) (2020), 612--640.
		
		\bibitem{OT91}  Ogawa, T.,  Tsutsumi, Y.: Blow-up of $H^1$ solution for the nonlinear Schr\"odinger equation. J. Differ. Equ. \textbf{92}(2) (1991), 317--330.
		
		\bibitem{S26}  Schr\"odinger, E.: Quantisierung als eigenwertproblem. Ann. Phys. \textbf{384}(4) (1926), 273--376.
		
		\bibitem{OS012}  Oh, T.,  Sulem, C.: On the one-dimensional cubic nonlinear Schr\"odinger equation below $L^2$. Kyoto J. Math.
		\textbf{52}(1) (2012), 99--115. 
		
		\bibitem{SS97} Sulem, C.,  Sulem, P-L.: Focusing nonlinear schr\"odinger equation and wave-packet collapse. Nonlinear Anal. \textbf{30}(2) (1997),  833--844.
		
		\bibitem{BS003}  Buslaev, V-S.,  Sulem, C.:  On asymptotic stability of solitary waves for nonlinear Schr\"odinger equations. Ann. I. H. Poincar\'{e} - AN. \textbf{20}(3) (2003), 419--475.
		
		\bibitem{TAX020}  Tuan, N.H.,  Au, V.V., Xu, R.: Semilinear Caputo time-fractional pseudo-parabolic equations. Commun. Pure Appl. Anal. \textbf{20}(2) (2021), 583--621.
		
		\bibitem{VB011} Vega, L., Banica, V.: Scattering for $1$D cubic NLS and singular vortex dynamics. J. Eur. Math. Soc. \textbf{14}(1) (2011), 209--253.
		
		\bibitem{BV020-1} Banica, V., Vega, L.: Evolution of polygonal lines by the binormal flow. Ann. PDE \textbf{6}(1) (2020), 1--53.
		
		\bibitem{BV020-2} Banica, V., Vega, L.: On the energy of critical solutions of the binormal flow. Commun. Partial Differ. Equa. \textbf{45}(7) (2020), 820--845.
		
		\bibitem{KPPV020} Kenig, C.E., Pilod, D., Ponce, G., Vega, L.: On the unique continuation of solutions to non-local non-linear dispersive equations, Commun. Partial Differ. Equ. \textbf{45}(8) (2020), 872--886.
		
		\bibitem{EKPV010} Escauriaza, L., Kenig, C.E., Ponce, G., Vega, L.: The sharp Hardy uncertainty principle for Schr\"odinger evolutions. Duke Math. J. \textbf{155}(1) (2010), 163--187.
		
		\bibitem{AV019} Agirre, M., Vega, L.: Some lower bounds for solutions of Schr\"odinger evolutions. SIAM J. Math Anal. \textbf{51}(4) (2019), 3324--3336.
		
		\bibitem{GV013} Gutierrez, S., Vega, L.: On the stability of self-similar solutions of $1$D cubic Schr\"odinger equations. Math. Ann. \textbf{356}(1) (2013), 259--300.
		
		\bibitem{AMV013} Alejo, M., Munoz, C., Vega, L.: The Gardner equation and the $L^2$-stability of the $N$-soliton solution of the Korteweg-de Vries equation. Trans. Amer. Math. Soc. \textbf{365}(1) (2013), 195--212.
		
		\bibitem{Xu017}  Xu, H.: Unbounded Sobolev trajectories and modified scattering theory for a wave guide nonlinear Schr\"odinger equation. Math. Z. \textbf{286} (2017), 443--489.
		
		\bibitem{W008}  Wang, Y.: Global existence and blow up of solutions for the inhomogeneous nonlinear Schr\"odinger equation in $\mathbf R^2$. J. Math. Anal. Appl. \textbf{338}(2) (2008), 1008--1019.

\bibitem{ZCL021}  Zhao, C.,  Caraballo, T.,  \L{}ukaszewicz, G.: Statistical solution and Liouville type theorem for the Klein-Gordon-Schr\"odinger equations.  J. Differ. Equ.
\textbf{281} (2021), 1--32.		
\end{thebibliography}
\end{document}